\documentclass[journal,12pt,onecolumn]{IEEEtran}
\usepackage{amsmath,amsfonts}
\usepackage{algorithmic}
\usepackage{algorithm}
\usepackage{array}
\usepackage[caption=false,font=normalsize,labelfont=sf,textfont=sf]{subfig}
\usepackage{textcomp}
\usepackage{stfloats}
\usepackage{url}
\usepackage{verbatim}
\usepackage{graphicx}
\usepackage{cite}
\usepackage{amssymb,epsfig,psfrag,tikz}
\usepackage{amsthm}
\usepackage{graphicx}
\usepackage{longtable}
\usepackage{array}
\usepackage[table]{xcolor}
\usepackage{booktabs}
\usepackage{color}
\usepackage{url}
\usepackage[colorlinks]{hyperref}
\hypersetup{
colorlinks=true,       
linkcolor=blue,          
citecolor=blue,        
filecolor=blue,      
urlcolor=blue 
}
\newtheorem{assumption}{Assumption}

\newtheorem{definition}{Definition}
\newtheorem{proposition}{Proposition}
\newtheorem{remark}{Remark}
\newtheorem{example}{Example}
\newtheorem{lemma}{Lemma}
\newtheorem{theorem}{Theorem}
\newcommand{\LCON}{L_{{\sf CO}}^{N}}
\newcommand{\LEXN}{L_{{\sf EX}}^{N}}
\newcommand{\LCOSN}{L_{{\sf CO,SYM}}^{N}}
\newcommand{\LCO}{L_{{{\sf CO}}}}
\newcommand{\LEX}{L_{{{\sf EX}}}}
\newcommand{\LCOS}{L_{{{\sf CO, SYM}}}}

\newcommand{\LPR}{L_{{{\sf PR}}}}
\newcommand{\LPRS}{L_{{{\sf PR, SYM}}}}
\newcommand{\LPRSN}{L_{{{\sf PR, SYM}}}^{N}}

\definecolor{darkred}{rgb}{0.7006, 0.03725, 0.03667}

\DeclareMathOperator*{\argminA}{arg\,min}

\allowdisplaybreaks

\begin{document}

\title{Decentralized Detection with Many Sensors: Optimality of Exchangeable and Identical Encoding Policies}

\author{Sina Sanjari$^{1}$, Naci Saldi$^{2}$, Sinan Gezici$^{3}$, Serdar Y\"uksel$^{4}$
\thanks{$*\,$The research of S. Sanjari and S. Yüksel is supported in part by the Natural Sciences and Engineering Research Council (NSERC) of Canada. S. Gezici acknowledges the support of the Scientific and Technological Research Council of Türkiye (TUBITAK) under the BIDEB 2219 Program.}
\thanks{$^{1}\,$Sina Sanjari is with the Department of Mathematics and Computer Science, Royal Military College of Canada, Canada;
        {\tt\small sanjari@rmc.ca}}%
\thanks{$^{2}\,$Naci Saldi is with the Department of Mathematics, Bilkent University, Ankara, Turkey;
        {\tt\small naci.saldi@bilkent.edu.tr}}%
\thanks{$^{3}\,$Sinan Gezici is with the Department of Electrical and Electronics Engineering, Bilkent University, Ankara, Turkey;
        {\tt\small gezici@ee.bilkent.edu.tr}}
\thanks{$^{4}\,$Serdar Y\"uksel is with the Department of Mathematics and Statistics, Queen's University, Canada;
        {\tt\small yuksel@queensu.ca}}
}


\maketitle

\begin{abstract}
We study a class of binary detection problems involving a single fusion center and a large or countably infinite number of sensors. Each sensor acts under a decentralized information structure, accessing only a local noisy observation related to the hypothesis. Based on this observation, sensors select policies to transmit a quantized signal through their actions to the fusion center, which makes the final decision using only these actions. This paper makes the following contributions: i) In the finitely many sensor setting, we provide a formal proof that an optimal encoding policy exists, and such an optimal policy is independent, deterministic, and of threshold type for the sensors and maximum \emph{a posteriori} probability type for the fusion center; ii) For the finitely many sensor setting, we further show that an optimal encoding policy exhibits an exchangeability (permutation invariance) property;  iii) We establish that an optimal encoding policy exists that is symmetric (identical) and independent across sensors in the infinitely many sensor setting under the error exponent cost; iv) Finally, we show that a symmetric optimal policy for the infinite population regime with the error exponent cost is approximately optimal for the large but finite sensor regime under the same cost criterion. We anticipate that the mathematical program used in the paper, building on \cite{sanjari2020optimality}, will find applications in several other massive communications applications.
\end{abstract}


\section{Introduction}
In decentralized detection problems, a collection of sensors gathers local observations from the environment and transmits quantized versions of these observations to a fusion center. The fusion center then decides on one of the hypotheses about the state of the environment based on the quantized observations coming from the sensors (and on additional observations, if any) \cite{tsitsiklis1989decentralizedDet,Varshney1997}. Numerous studies have been conducted on decentralized detection, and both its theoretical and practical aspects have been examined in various studies
\cite{tsitsiklis1989decentralizedDet,Varshney1997,ViswanathanVarshney1997,BlumKassamPoor1997,Tenney1981}.

The primary objective in a decentralized detection problem is to design the decision rules at both the sensors and the fusion center in an optimal manner to minimize a given cost function (such as the average Bayes risk). Once the decision rules at the sensors are determined (fixed), obtaining the optimal decision rule at the fusion center is not challenging, as it reduces to the classical design of optimal decision rules \cite{PoorBook,Chair1986}. For instance, under the criterion of minimizing the average probability of error, using the maximum \emph{a posteriori} probability (MAP) decision rule at the fusion center is optimal for any given sensor outputs \cite{tsitsiklis1989decentralizedDet,Drakopoulos1991,Kam1992,Chair1986}. However, when the decision rule at the fusion center is fixed, determining the optimal decision rules at the sensors generally becomes a highly complex and challenging problem. This complexity arises due to two main reasons. First, although sensors employ separate decision rules, their observations may be statistically dependent, severely complicating the optimal design of sensor decision rules. Second, even for independent sensor observations, obtaining the globally optimal solution for the joint design of all sensor decision rules corresponds to a complex optimization problem in general. Therefore, in the literature, most studies adopt the person-by-person optimality criterion and focus on scenarios in which sensor measurements are independent under each hypothesis \cite{tsitsiklis1989decentralizedDet,Varshney1997}. In this setting, computing the optimal decision rules at the sensors becomes more manageable as they can be represented by threshold tests \cite{ViswanathanVarshney1997,Warren1999,varshney2020decentralized}. Also, extensions to distributed detection problems with noisy communication channels between sensors and the fusion center have been reported \cite{Chen2009,DecentCSI}.

When the sensor observations are dependent, no general and computationally feasible solution exists in the literature for obtaining the optimal decision rules at the sensors. In particular, even for the binary hypothesis testing problem with the person-by-person optimality criterion, there does not exist a finite parametrization of the candidate decision rules as threshold strategies may not be optimal \cite{Yan2001,Tang1992}. Hence, the problem is inherently intractable in general \cite[Sec.~2.5]{tsitsiklis1989decentralizedDet}, \cite{Tsitsiklis1985}. However, there exist few results in the literature for certain special cases \cite{Chen1995,Willett2000,DDcorre,ChenVarshney2012}. For instance, in \cite{ChenVarshney2012}, a hierarchical conditional independence model is assumed, leading to results similar to those in the independent case. However, this assumption imposes a significant constraint on the problem structure. Also, \cite{Chen1995} and \cite{Willett2000} consider a specific problem of binary detection with two sensors that observe a shift in the mean of correlated Gaussian random variables. As an alternative approach, asymptotic analyses can be performed by considering the error exponent as the detection performance metric  \cite{Han1994,Han1998,Ahlswede1986}.

As noted, the current literature primarily focuses on person-by-person optimality in which the decision rule of each sensor is optimized when all the other decision rules are fixed \cite{varshney2020decentralized}. Although this approach simplifies the design of sensor decision rules, it may not lead to a locally or globally optimal solution in general \cite{tsitsiklis1989decentralizedDet}. There exists a limited number of studies that focus on global optimality instead of person-by-person optimality for decentralized detection problems \cite{Tenney1981}, \cite[section~3]{varshney2020decentralized}, \cite{DDglobal1987,WarrenWillett,tsitsiklis1988decentralized}. For example, minimizing a minimum Bayes risk expression over the sensor threshold parameter is performed in \cite[Sec.~3]{varshney2020decentralized} by considering identical sensor decision rules and the optimal Bayes rule at the fusion center. Also, the authors of \cite{DDglobal1987} focus on a decentralized binary hypothesis testing problem where local sensors send their binary (one bit) decisions to a fusion center, and optimize the decision rules at both the local sensors and the fusion center based on the approaches in \cite{Tenney1981} and \cite{Chair1986}. In particular, the optimal sensor decision rules are specified as likelihood ratio tests (LRTs) in \cite{Tenney1981} by considering a fixed decision rule at the fusion center. The study in  \cite{WarrenWillett} considers multi-bit sensor decision rules and proves that the optimal rule at a sensor is a partition determined by the likelihood ratio at the sensor input according to both the Neyman-Pearson and Bayes criteria. It is also shown that sensor decision rules based on likelihood ratio partitions maximize the Ali-Silvey distance between the distributions under the two hypotheses at the fusion center
input \cite{WarrenWillett}. 

{We also note that recent work has considered decentralized detection under wireless constraints. 
\cite{Ciuonzo2025ScoreBased} investigates the detection of an unknown deterministic signal in a wireless sensor network where sensors send quantized measurements over Rayleigh fading channels, and develops score tests that jointly perform decoding and fusion with low complexity. 
\cite{Mao2025Hybrid} studies signal detection under additive and multiplicative Gaussian noise with both quantized and full-precision sensor observations, deriving a hybrid detector together with quantizer-threshold optimization and bandwidth allocation. \cite{CiuonzoRossi2018Quantizer} addresses target detection with unknown emitted power and position, using locally optimum fusion rules and sensor-threshold optimization under binary channels. These works focus on fusion sensor design and resource allocation for wireless sensor networks, whereas our focus is on existence, structural optimality, exchangeability, and symmetric optimality in many-sensor decentralized detection.}

\begin{table}[p]
\centering
\caption{Summary and setting comparison of related works.}
\label{tab:related-works-combined}

\begingroup
\scriptsize
\setlength{\tabcolsep}{2pt}
\renewcommand{\arraystretch}{0.85}



\resizebox{\textwidth}{!}{%
\begin{tabular}{
@{}
>{\raggedright\arraybackslash}p{0.05\textwidth}
>{\raggedright\arraybackslash}p{0.4\textwidth}
>{\raggedright\arraybackslash}p{0.55\textwidth}
@{}
}

\toprule
\textbf{Related work} & \textbf{Setting} & \textbf{Most Relevant Contribution} \\
\midrule

\cite{tsitsiklis1989decentralizedDet}
&
Decentralized detection under independent and dependent observation models.
&
Foundational formulation; threshold-structure results under conditional independence; discussion of dependent-observation difficulties, computational issues, and asymptotic results.\\

\addlinespace

\cite{Varshney1997}
&
Decentralized detection under various sensor topologies.
&
Reviews decentralized detection under various sensor topologies and discusses structural results, local decision-rule design, fusion-rule design, and globally optimal formulations.
\\

\addlinespace

\cite{Tenney1981}
&
Conditionally independent binary decentralized hypothesis testing problem with two sensors and a fixed fusion rule.
&
Shows that the solution to the problem is deterministic and involves a likelihood ratio test for each detector. Also, provides the computations of the thresholds for these tests that are generally coupled but may be only locally optimal.\\

\addlinespace

\cite{Chair1986}
&
Fusion-center design for fixed local sensors' decisions under conditionally independent observations.
&
Derives optimal fusion rules when the individual detector decision rules are known, and shows that under this policy, individual decisions are
weighted according to their reliability, weighted as a function of the probability of miss and false alarm of the individual detector. \\

\addlinespace

\cite{Drakopoulos1991}
&
Decentralized binary hypothesis testing with correlated local sensor decisions represented through a finite set of correlation coefficients.
&
Derives the optimum decision fusion rule in the Neyman-Pearson sense, shown to be a function of correlation coefficients, and analyzes the performance based on the correlation between sensors' decisions. \\

\addlinespace

\cite{Kam1992}
&
Fusion-center design for fixed correlated local binary decisions.
&
Develops optimal data-fusion rules for correlated local binary decisions in terms of correlation coefficients of different orders. \\

\addlinespace

\cite{ViswanathanVarshney1997}
&
Decentralized detection under parallel and serial sensor topologies based on Neyman-Pearson criterion and Bayesian risk with conditionally independent sensor observations.
&
Basic results on distributed detection are reviewed, including optimality of LRTs and derivations of threshold- or likelihood-ratio-based local sensor rules. \\

\addlinespace

\cite{Warren1999}
&
Decentralized detection for Neyman-Pearson, Bayes, Ali-Silvey distance, and mutual information criteria with conditionally independent sensor observations.
&
Shows that optimal sensors' policies use a likelihood-ratio partition of the sensor observations. Also, shows for conditionally i.i.d. observations, the policies used for the sensors are identical under a subclass of the Ali-Silvey distances or local mutual information criterion.\\

\addlinespace

\cite{Chen2009}
&
Decentralized decision rules under conditionally independent observations, where sensors' decisions are transmitted over dependent nonideal channels.
&
Under Bayesian and Neyman-Pearson criteria, optimal detection is achieved using likelihood-ratio quantizers at local sensors, for various channel characteristics, including correlated channels between the sensors and the fusion center. \\

\addlinespace

\cite{DecentCSI}
&
Local decision-rule design for decentralized detection over fading channels with either global channel-state information or channel-fading statistics.
&
Derives the necessary conditions for optimal local sensor decision rules, and shows the optimality of the likelihood ratio test. \\

\addlinespace

\cite{DDcorre}
&
Detection under a known signal with correlation between sensors, and also an unknown signal in white noise.
&
Restricting to detection policies based on local LRT, develops LRT-type rules and determines the performance. \\

\addlinespace

\cite{Tang1992}
&
$M$-ary hypothesis testing 
problem with correlated observations.
&
Derives the necessary conditions 
for the person-by-person optimal policies and uses Gauss--Seidel iterations to obtain locally optimal rules. \\

\addlinespace

\cite{Yan2001}
&
Neyman--Pearson decentralized detection with correlated observations.
&
Provides necessary conditions for Neyman-Pearson optimal detectors for a parallel sensor topology, and also provide a procedure for finding the Neyman-Pearson optimal policy.  \\

\addlinespace

\cite{Tsitsiklis1985}
&
Computational complexity of decentralized detection (due to non-classical information).
&
Establishes complexity and intractability results, especially when the assumption of conditional independence is removed. \\

\addlinespace

\cite{ChenVarshney2012}
&
Decentralized detection with a hierarchical conditional-independence model.
&
By introducing a hidden variable, classes of distributed detection problems with dependent observations whose optimal sensor decision rules are uncoupled in their form are identified. \\

\addlinespace

\cite{Chen1995}
&
Two-sensor correlated Gaussian detection.
&
Finds a sufficient condition on the noise mean and covariance under which the optimal binary quantizers are threshold-type functions of the marginal likelihood ratio. \\

\addlinespace

\cite{Willett2000}
&
Detection of a mean shift in a set of correlated Gaussian observations with two sensors. 
&
Characterizes regimes in which threshold-type rules are effective or fail. \\

\addlinespace

\cite{DDglobal1987}
&
Joint local sensor and fusion-rule optimization under conditional independence.
&
Derives coupled optimality equations for global design and compares the optimal performance
with two suboptimal ones. \\

\addlinespace

\cite{WarrenWillett}
&
Local decisions under Neyman-Pearson and Bayes criteria with conditionally independent observations.
&
Shows that an optimal sensor policy under both the Neyman-Pearson and the Bayes criteria is a partition determined by the likelihood ratio at the sensor input. \\

\addlinespace

\cite{tsitsiklis1988decentralized}
&
Large-sensor asymptotics for $M$-ary hypothesis testing with conditionally i.i.d. observations and $D$-valued messages.
&
Shows asymptotic optimality of partitioning sensors into $M(M-1)/2$ groups for $M$-ary case (shows that identical LRTs are asymptotically optimal when $M=D=2$). \\

\addlinespace

\cite{CiuonzoRossi2018Quantizer}
&
Detection over error-prone binary symmetric channels between sensors and the fusion center with a wireless sensor network.

&
Develops locally optimum detector fusion and channel-aware threshold optimization. \\

\addlinespace

\cite{Ciuonzo2025ScoreBased}
&
Decentralized detection of an unknown deterministic signal in a wireless sensor network.
&
Proposes score tests as low-complexity alternatives to the generalized likelihood ratio
test.  \\

\addlinespace

\cite{Mao2025Hybrid}
&
Distributed detection with hybrid quantized and full-precision observations over multiplicative fading channels.
&
Optimizes quantization thresholds and bandwidth allocation for hybrid quantized/full-precision detection.\\

\bottomrule

\end{tabular}
}

\endgroup

\end{table}

In a seminal work by Tsitsiklis \cite{tsitsiklis1988decentralized}, an asymptotic analysis is presented, where sensors with independent and
identically distributed observations generate $D$-valued messages and send them to a fusion center in an $M$-ary hypothesis testing framework. It is shown that as the number of sensors goes to infinity, it is asymptotically
optimal to divide the sensors into $M(M-1)/2$ groups, with all sensors in each group employing the same decision rule. For $M=D=2$, the asymptotic solution simplifies to having all the sensors perform an identical LRT  \cite{tsitsiklis1988decentralized}. {A summary of the literature review presented above is provided in Table~\ref{tab:related-works-combined} together with setting comparisons.}

Our aim in this work is to focus on the global optimization of sensor decision rules in a decentralized binary detection problem and to establish existence and structure results for optimal decision rules, considering both finite and infinite numbers of sensors. To achieve this goal, we utilize recent results for large stochastic teams in \cite{sanjari2020optimality}. In stochastic team decision problems, a finite number of decision-makers (agents) act collectively to minimize a common cost function, but they do not necessarily share their information; see, e.g.,  \cite{YukselBasar24,Radner1962}. In \cite{sanjari2020optimality}, optimality of independently randomized symmetric policies has been established for a class of exchangeable stochastic teams. We view the decentralized detection problem as a stochastic team problem involving the sensors and the fusion center, building also on the mathematical machinery developed in \cite{sanjari2020optimality} while accounting for the unique structure of the detection problem.

In this paper, we derive several results on the existence and structure of optimal policies for a decentralized detection problem, assuming that the observations are independent and identically distributed (i.i.d.) among sensors, given the hypothesis. We also provide asymptotic analysis as the number of sensors approaches infinity. The main contributions of the paper can be summarized as follows:
\begin{enumerate}
    \item In Theorem \ref{the:thre}, we establish the existence of an optimal policy and that such an optimal policy is deterministic, and of threshold type for sensors in the finite population setting. Optimality of threshold policies (assuming existence) was established in the seminal work by Tsitsiklis (see \cite[Proposition 2.4]{tsitsiklis1989decentralizedDet}); here, we provide an existence proof in Theorem \ref{the:thre}. Key to our analysis in this case is a reformulation of the problem and the approach via an optimal quantization formulation, considering optimality of threshold policies established in \cite{tsitsiklis1989decentralizedDet}.
    \item For the finitely many sensor setting, in Theorem \ref{the:exi-N}, we show that an optimal encoding policy exists and that such an optimal policy is exchangeable among sensors. In this formulation, the fusion center is allowed to have access to the possible randomness of the encoders' policies.
    \item In Theorem \ref{the:Exist-inf}, we demonstrate that an optimal solution exists and that such an optimal policy is symmetric (identical) and independent among sensors in the infinite population regime under the error exponent cost. 
    \item In Proposition \ref{the:thre-inf}, we establish that a symmetric optimal policy for the infinite population regime with the error exponent cost is approximately optimal for the problem with a large number of sensors with the same cost.
\end{enumerate}

The remainder of the paper is organized as follows. In Section~\ref{sec:Binary}, we formally introduce the problem formulation and
summarize the technical contributions. 
In Section \ref{sec:exitsN}, Section \ref{sec:Randpolicies}, and Section \ref{se:ExitExN}, we introduce various notions of policies and establish the existence and structural results for an optimal policy for a finite population setting. In Section \ref{sec:Inf}, we present the results for the infinite population setting.

\section{Binary Detection Problem with Uniform Cost Assignment}\label{sec:Binary}

We study a class of binary detection problems in the Bayesian framework with uniform cost assignment (UCA), where the possible hypotheses are denoted by $H_1$ and $H_2$. The true hypothesis $H^{\star}$ can be regarded as a random variable equal to $H_j$ with probability $P(H_j)$ for $j=1,2$. We focus on problems with a finite number of sensors as well as a countably infinite number of sensors.

 In the finite population setting, there exist $N$ peripheral sensors in the environment, which make observations from the environment with $y^i$ denoting the observation (measurement) of sensor $i$, taking values in an appropriate dimensional Euclidean space $\mathbb{Y}^i$ for $i\in\{1,\ldots,N\}$. Each sensor quantizes its observation via a mapping $\gamma^i$ that maps the observation $y^i$ into an action $u^i$. In particular, for a given realization $y^i$, an action of the sensor $i$ is expressed as $u^i=\gamma^i(y^i)\in \mathbb{U}^i=\{1,\ldots,|\mathbb{U}^i|\}$ for $i\in\{1,\ldots,N\}$, where $\gamma^i:\mathbb{Y}^i\rightarrow \mathbb{U}^i$ is a policy or \textit{decision rule of the sensor $i$}. Each sensor's action is observed by the fusion center (denoted by agent $0$), which makes the final decision about the true hypothesis based on the sensors' actions. Therefore, the final decision is stated as
$u^0=\gamma^0(u^1,\ldots,u^N)$, where $\gamma^0:\prod_{i=1}^{N}\mathbb{U}^i \rightarrow \{1,2\}$ is the policy (or decision rule) of the fusion center.  We use $\Gamma^i$ to represent the set of all possible policies for agent $i$ for $i=0,1,\ldots,N$. The collection of policies is denoted by $(\gamma^0,\gamma^1,\ldots,\gamma^N)$. The set of all policies is expressed as $\Gamma^0\times\Gamma^1\times\cdots\times\Gamma^N$.  The information structure of the detection problem with $N$ sensors is depicted in Fig. \ref{figure:0}.

\begin{figure}[h]
    \centering
     \scalebox{.8}{
    \begin{tikzpicture}[node distance=5cm, auto]
        
        \node (E) [draw, rectangle, minimum width=1cm, minimum height=1cm, yshift=0.75cm] {$H$};
        \node (A) [draw, circle, minimum width=1cm, minimum height=2cm, right of=E, xshift=-1cm, yshift=0cm,text width=1.8cm] {Observation\\ \:\:\:channel};
        
        \node (B1) [draw, rectangle, minimum width=2cm, minimum height=1cm, right of=A, yshift=2cm] {Sensor $1$};
        \node (B2) [draw, rectangle, minimum width=2cm, minimum height=1cm, right of=A, yshift=0.5cm] {Sensor $2$};
        \node (Bdots) [draw=none, right of=A, yshift=-0.5cm] {$\vdots$};
        \node (BN) [draw, rectangle, minimum width=2cm, minimum height=1cm, right of=A, yshift=-1.5cm] {Sensor $N$};
        \node (C) [draw, rectangle, minimum width=2cm, minimum height=1cm, right of=Bdots, yshift=0.5cm] {Fusion Center};
        \node (D) [text width=1cm] at (17.5,0.76) {};
        \node[text width=1cm] at (7.8,2.6) {$y^1$};
        \node[text width=1cm] at (7.8,1.5) {$y^2$};
        \node[text width=1cm] at (7.8,0) {$y^N$};
        \node[text width=1cm] at (11,2.6) {$u^1$};
        \node[text width=1cm] at (11,1.5) {$u^2$};
        \node[text width=1cm] at (11,0) {$u^N$};
        \node[text width=1cm] at (16.4,1.2) {$u^0$};
        
        \draw[->] (E) -- (A);
        \draw[->] (A) -- (B1);
        \draw[->] (A) -- (B2);
        \draw[->] (A) -- (BN);
        \draw[->] (B1) -- (C);
        \draw[->] (B2) -- (C);
        \draw[->] (BN) -- (C);
        \draw[->] (C) -- (D);
    \end{tikzpicture}
    }
    \caption{The information structure of the decentralized detection problem with $N$ sensors. Each sensor observes $y^i$ and sends information $u^i$ to the fusion center.}\label{figure:0}
\end{figure}

We assume that the action and observation spaces are identical among sensors, denoted by $\mathbb{U}$ and $\mathbb{Y}$, respectively. We consider the cost function $c$ that depends directly only on the fusion center's decision, $u^0$. In particular, we focus on the following cost:
\begin{gather}\label{eq:probErr}
c(H_j,u^0)=
\begin{cases}
1,&{\rm{if}}~u^0\ne j\\
0,&{\rm{if}}~u^0= j
\end{cases}.
\end{gather}

The above detection formulation is an instance of a stochastic team problem with the following expected cost:
\begin{align}\label{eq:costStoc}
    J^{N}(\gamma^{0:N})=\mathbb{E}^{\gamma^{0:N}}\left[c(H,u^{0})\right].
\end{align}
We call the problem involving the cost in \eqref{eq:costStoc} as $\mathcal{P}^{N}$. Hence, we have a dynamic stochastic team problem with a non-classical information structure; see \cite[Section 3.2]{YukselBasar24}. The information structure is non-classical because the fusion center observes the actions of the sensors but not their observations. Our approach is to view the decentralized detection problem, with the fusion rule set as the optimal rule, and thus remove it from the system as a static team problem. We characterize the fusion center's optimal decoding rule for any sensor's encoding policy. 

For the infinite population problem, following \cite[Proposition 5.1]{tsitsiklis1989decentralizedDet}, the expected cost (which is the probability of error) goes to zero exponentially as the number of sensors goes to infinity. Hence, we instead study the error exponent as our new cost for the infinite population problem, i.e., the error exponent cost in the finite population for any policy $\gamma^{0:N}$ is given by
\begin{align}\label{eq:JEE}
 J^N_{\sf EE}(\gamma^{0:N})&= \frac{\log \left(J^N (\gamma^{0:N})\right)}{N}
\end{align}
and the error exponent cost of a policy $\gamma^{0:\infty}:=(\gamma^0, \gamma^1, \ldots)$ 
\begin{align}\label{eq:JEEinf}
    J_{\sf EE}(\gamma^{0:\infty})&=\limsup_{N\to \infty} \frac{\log \left(J^N(\gamma^{0:N})\right)}{N}.
\end{align}
We call the problems involving the costs in \eqref{eq:JEE} and \eqref{eq:JEEinf} by $\mathcal{P}^{N}_{\sf E}$ and $\mathcal{P}^{\infty}_{\sf E}$, respectively.

Next, we introduce our main assumptions.

\begin{assumption}\label{Assump:density}
Suppose that
\begin{itemize}
    \item [(i)] $y^1, \ldots, y^N$ are i.i.d., conditioned on $H$.
\item [(ii)] $\mathbb{U}^i=\mathbb{U}$ for $i=1, \ldots, N$, where $\mathbb{U}$ is a finite set.
\item[(iii)] The conditional density $f_{y^i|H}$ of $y^i$ given $H$ for every $i\in\{1,\ldots,N\}$ exists,
$f_{y^i|H}(y^i|H_1)>0$ for all $y^i\in\mathbb{Y}$, and the likelihood ratio
\begin{equation}\label{eq:L^i}
 l^i = L(y^i) := \frac{f_{y^i|H}(y^i|H_2)}{f_{y^i|H}(y^i|H_1)} \qquad a.s. 
\end{equation}
is continuously differentiable in $y^i$ with a gradient that does not vanish, i.e., $\nabla L(y^i)\neq 0$ for all
$y^i\in\mathbb{Y}$.
\end{itemize}
\end{assumption}

Assumptions \ref{Assump:density}-(i) and \ref{Assump:density}-(iii) allow us to probabilistically view the observations of sensors as independent of each other and $H$; see \cite{alhakeem1995unified,varshney2020decentralized}. Additionally, under Assumption~\ref{Assump:density}--(iii), the random variable $l^i: L (y^i)$ admits a density using the inverse function theorem of vector calculus \cite[Section 1.11]{hajek2015random}.

Our main goal is to establish the existence and structural results for an optimal solution for the $N$-sensor setting and its limiting problem with an infinite number of sensors. To this end, we allow joint randomization among the sensors' policies and expand the information of the fusion center to include the realization of joint sensors' policies (common randomness). 
Thus, we allow the fusion center to have access to the (\emph{known common randomization}) information set $I^{\sf KCR}_N$ with
\begin{align}\label{eq:IS-KR}
    I^{\sf KCR}_N:=\{u^{1:N}, \gamma^{1:N}\},
\end{align}
under which the common randomness is known to the fusion center. For our analysis, we work with another information set for the fusion center, which is amenable to our convergence argument as $N\to \infty$. 

{A summary of our notation is included in the table below.}

\begingroup
\scriptsize
\setlength{\tabcolsep}{1pt}
\renewcommand{\arraystretch}{1.2}

\begin{longtable}{
>{\raggedright\arraybackslash}p{0.1\textwidth}
>{\raggedright\arraybackslash}p{0.38\textwidth}
>{\raggedright\arraybackslash}p{0.1\textwidth}
>{\raggedright\arraybackslash}p{0.38\textwidth}
}

\caption{Summary of notation.}
\label{tab:notation-summary}
\\

\toprule
\textbf{Notation} & \textbf{Description}
& \textbf{Notation} & \textbf{Description} \\
\midrule
\endfirsthead

\toprule
\textbf{Notation} & \textbf{Description}
& \textbf{Notation} & \textbf{Description} \\
\midrule
\endhead

\multicolumn{4}{l}{\textbf{Basic detection model}} \\
\midrule

$N$
& Number of sensors in the finite-population problem.
&
$H_1,H_2$
& The two hypotheses in the binary detection problem. \\

$y^i$, $\mathbb{Y}$
& Observation of sensor $i$ and its observation space.
&
$u^i$, $\mathbb{U}$
& Action of sensor $i$ and the finite sensor-action alphabet. \\

$u^0$, $\mathbb{U}^0$
& Fusion-center action and its action alphabet.
&
$c$
& Detection cost. \\

\addlinespace[2pt]
\multicolumn{4}{l}{\textbf{Policies and policy spaces}} \\
\midrule

$\gamma^i$, $\Gamma^i$
& Encoding policy of sensor $i$ and its admissible policy space.
&
$\gamma^0$, $\Gamma^0$
& Decoding policy of the fusion center and its policy space. \\

$\gamma^R$, $\Gamma^R$
& Representative-sensor policy and policy space.
&
$\Gamma^N_{\sf TS}$
& Class of deterministic threshold-type sensor policies. \\

\addlinespace[2pt]
\multicolumn{4}{l}{\textbf{Cost criteria and problem formulations}} \\
\midrule

$\mathcal{P}^N$
& Finite-population problem under probability-of-error cost.
&
$\mathcal{P}^N_E$
& Finite-population problem under error-exponent cost. \\

$\mathcal{P}^\infty_E$
& Infinite-population problem under error-exponent cost.
&
$J^N$
& Expected detection cost in the finite-sensor problem. \\

$J^N_{\sf EE}$
& Error-exponent cost for the finite-sensor problem.
&
$J_{\sf EE}$
& Error-exponent cost for the infinite-population problem. \\

\addlinespace[2pt]
\multicolumn{4}{l}{\textbf{Likelihood statistics and information variables}} \\
\midrule

$l^i$
& Likelihood ratio associated with sensor $i$'s observation.
&
$\Delta_N$
& Log-likelihood-ratio statistic used by the fusion center. \\

$I^{\sf KCR}_N$
& Fusion-center information when common randomization is known.
&
$I^{\sf FS}_N$
& Fusion-center information in the finite-sensor randomized formulation. \\

$\mu^N$
& Empirical measure of the realized sensor policies.
&
&
\\

\addlinespace[2pt]
\multicolumn{4}{l}{\textbf{Randomized policy classes}} \\
\midrule

$L^N$
& Randomized joint policies for finitely many sensors.
&
$L$
& Randomized policies for countably many sensors. \\

$L^N_{\sf CO}$
& Finite-sensor policies with common and independent randomization.
&
$L_{\sf CO,SYM}$
& Symmetric policies with common and independent randomization. \\

$L^N_{\sf EX}$
& Exchangeable randomized finite-sensor policies.
&
$L_{\sf EX}$
& Infinitely exchangeable randomized policies. \\

$L^N_{\sf PR}$
& Independently randomized finite-sensor policies.
&
$L^N_{\sf PR,SYM}$
& Independently randomized and symmetric finite-sensor policies. \\

$L_{\sf PR,SYM}$
& Symmetric independently randomized policies for infinitely many sensors.
&
&
\\

\bottomrule

\end{longtable}

\endgroup
\color{black}

Before getting into the details of the theoretical results, we summarize the main contributions in the next section.

\subsection{Summary and Discussion of Main Results and Technical Approach}\label{sec:Summary}

 \begin{figure}[ht]
\centering
    \begin{tikzpicture}[
     very thin]
    \filldraw[color=blue!60, fill=blue!15, very thick] (-3,0.8) rectangle (2.3,2.5);
    
    \filldraw[color=blue!60, fill=blue!15, very thick] (4.5,0.8) rectangle (10.3,2.5);

    \draw[gray, thick] (3.5,1.1) -- (3.5,-4.3);
    \draw[gray] (-3.8,2.4) -- (-3.8,-4.3);
    \draw[gray] (11.2,2.4) -- (11.2,-4.3);

    \draw[gray] (-3.8,0.6) -- (11.2,0.6);

    \node[text width=5cm] at (0,1.6) 
    {Detection problem with $N$ exchangeable sensors};

    \node[black,text width=7cm] at (0,-0.5) 
    {\small{{\color{blue}{$\blacksquare$}}\: \sf Theorem \ref{the:thre}: Existence of an optimal policy that is of threshold type for sensors and MAP for the fusion center.}};
\node[black,text width=7cm] at (0,-2.9) 
    {\small{{\color{blue}{$\blacksquare$}}\: \sf Lemma \ref{lemma:1-EX} and Theorem \ref{the:exi-N}: Existence of an optimal policy that is exchangeable for sensors and MAP for the fusion center.}};

\draw[blue, thick, ->] (2.8,1.6) -- (4.2,1.6);
\node[black,text width=2cm] at (3.8,2) 
    {{\color{blue} $N$ $\rightarrow\infty$}};

    \node[text width=5cm] at (7.4,1.6) 
    {Detection problem with infinite exchangeable sensors};
    
    \node[black,text width=7cm] at (7.4,-0.7) 
    {\small{{\color{blue}{$\blacksquare$}}\: \sf Theorem \ref{the:Exist-inf}: Existence of an optimal policy that is symmetric and independent (possibly randomized) for sensors and MAP for the fusion center.}};
    \node[black,text width=7cm] at (7.4,-3) 
    {\small{{\color{blue}{$\blacksquare$}}\: \sf Proposition \ref{the:thre-inf}: Asymptotic optimality of a symmetric and independent optimal policy for the problem with a finite but large number of exchangeable sensors.}};
\end{tikzpicture}
\caption{Summary of main contributions.}
\label{fig:theorem-graph}
\end{figure}

The main contributions of the paper are summarized in Fig. \ref{fig:theorem-graph}. We first consider the detection problem with $N$ sensors under deterministic encoding policies. We define threshold policies for sensors in Definition \ref{def:threshold}. Then, via Lemmas \ref{lem:TS-optimal} and \ref{lemma:2-PI}  \cite{tsitsiklis1989decentralizedDet}, we establish that without any loss, the search for an optimal encoding policy of the sensors can be restricted to threshold policies and the search for an optimal decoding policy of the fusion center to MAP rules. Given this, we present an equivalent problem (see Fig. \ref{figure:2}) in which each sensor only utilizes the ratio term in \eqref{eq:L^i} and the fusion center utilizes the log likelihood ratio term. The optimality results in Lemmas \ref{lem:TS-optimal} and \ref{lemma:2-PI} do not guarantee the existence of an optimal solution. To guarantee the existence of an optimal solution, we view the set of threshold policies as quantizer maps which are then viewed as stochastic kernels with a fixed input marginal, which we then show to form a compact space under total variation building on \cite{yuksel2012optimization}; see Lemma \ref{lem:compactness of TS}. Leveraging this compactness, we establish the existence of an optimal threshold policy in Theorem \ref{the:thre} by proving that, for any MAP rule at the fusion center, the cost is continuous with respect to the sensors’ threshold policies under total variation.

To establish a connection to an infinite population under the error exponent cost, we allow joint randomization among sensors, where we introduce various sets of randomized policies (see Fig. \ref{Fig-policies}). To show that such a randomization is legitimate in the sense that it does not change the performance of the problem, we expand the information set of the fusion center (see information sets \eqref{eq:IS-KR} and \eqref{eq:IFS}). In Lemma \ref{Lem:info-optimal-N}, we show that under these equivalent information sets, the randomization in encoding policies does not change the optimal performance. We then, in Theorem \ref{the:exi-N}, show that among all randomized encoding policies, there exists an optimal solution for the finite sensor problem that is exchangeable among sensors (the joint distribution on policies is permutation-invariant). As it is demonstrated in Example \ref{example1}, an optimal solution might not be symmetric (identical) and independent. 

Utilizing a convergence property of exchangeable sequences in Theorem \ref{the:ald} by Aldous \cite[Proposition 7.20]{aldous2006ecole}, we show in Theorem \ref{the:Exist-inf} that a sequence of optimal exchangeable policies for the finite agent problem under error exponent cost converges in distribution (under a subsequence) to an infinitely exchangeable policy (which is a mixture of independent and symmetric policies by the de Finetti representation Theorem \ref{the:defin} \cite{aldous2006ecole}). This additionally establishes that the induced empirical measures of joint observations, actions, and policies converge weakly to a limit that is computable by the fusion center. This then establishes the optimality of the limiting policy utilizing a lower bound for the probability of error in  \cite{shannon1967lower}. Finally, using a convex analytic approach, in Theorem \ref{the:Exist-inf} for the infinite sensor problem, we establish that an optimal policy exists and is symmetric and independent for sensors and an MAP rule for the fusion center. This, in turn, leads to optimality of an equivalent problem with a single representative sensor and a fusion center problem (see \eqref{eq:rep-encod}). Additionally, this result leads to an interesting conclusion, demonstrating that a lack of knowledge about the randomization of sensors' policies does not incur any loss in performance for the detection problem asymptotically as the number of sensors approaches infinity. As another implication of our result, in Proposition \ref{the:thre-inf}, we show that a symmetric optimal solution of the representative problem is approximately optimal for the problem with a finite but large number of exchangeable sensors.

\section{Equivalent Formulation for the Finitely Many Sensor Problem}\label{sec:exitsN}

We first define a class of deterministic \emph{threshold} policies for sensors adopted from \cite{tsitsiklis1989decentralizedDet}.

\begin{definition}\label{def:threshold}
 
\begin{enumerate}
    \item A deterministic policy for sensor $i$, $\gamma^i$, is a monotone threshold policy if there exist an integer $2\leq m^i\leq |\mathbb{U}|$ and threshold constants $\{t_1^i,\ldots,t_{m^i-1}^i\}$  satisfying $0\leq t_1^i \leq t_2^i \leq \cdots \leq t_{m^i-1}^i \leq \infty$ such that 
\begin{gather}
\gamma^i(y^i)=d \qquad \text{if}\quad l^i\in B_d^i
\end{gather}
for $B_d^i=[t_{d-1}^i,t_d^i)$ with $d=2, \ldots,  m^i-1 $, $B_1^i=[0,t_1^i)$, and $B_{m^i}^i=[t_{m^i-1}^i,\infty]$, where
\begin{gather}\label{eq:li}
    l^i:=L(y^i):=\frac{f_{Y^i|H}(y^i|H_2)}{f_{Y^i|H}(y^i|H_1)}\qquad a.s.
\end{gather}

\item A deterministic policy for a sensor $i$,
 $\gamma^i$, is a threshold policy if there exists a permutation mapping $\tau:\mathbb{U}\to\mathbb{U}$  such that $\tau \circ \gamma^i$
is a monotone threshold policy as defined above.
\end{enumerate}
\end{definition}
{We note that a monotone threshold policy is a special case of a threshold policy in which the output label increases with the likelihood-ratio interval. Under a general threshold policy, the likelihood-ratio space is still partitioned into ordered threshold intervals; however, the labels assigned to these intervals need not be monotone in their original order, but are monotone only after applying a permutation $\tau$ of the output labels.} Denote the space of all deterministic threshold policies by $\Gamma_{\sf TS}^N$.  We note that allowing the number of bins $m^i$ to vary will be useful for establishing compactness properties, which will be discussed later in the paper under a suitable topology on the space of sensor policies. The optimality of threshold policies is stated next. The proof of the following lemma for completeness is included in Appendix \ref{APP:lem TS}.

\begin{lemma}\label{lem:TS-optimal}\cite{tsitsiklis1989decentralizedDet}
    Consider $\mathcal{P}^{N}$. Let Assumption \ref{Assump:density} hold. Then,
    \begin{align*}
    \inf_{\gamma^{1:N}\in \Gamma_N} J^N(\gamma^0,\gamma^{1:N})&=\inf_{\gamma^{1:N}\in \Gamma_{\sf TS}^N} J^N(\gamma^0,\gamma^{1:N})
\end{align*}
for any $\gamma^0\in \Gamma^0$.
\end{lemma}
We now provide the main reasoning behind the proof of Lemma \ref{lem:TS-optimal}. Suppose that we fix policies of sensors with index $k\in \{1, \ldots, i-1, i+1, \ldots, N\}$ to deterministic policies $\gamma^{k\star}$. The optimization problem  for the deviating sensor $i$ is given by
\begin{align*}
 \inf_{u^i\in \mathbb{U}} \left\{\sum_{j=1}^{2} P(H_j|y^i)\mathbb{E}^{\gamma^{-i\star}}\left[c(H_j,\gamma^{0}(u^{1:N}))\mid H_j\right]\right\},
\end{align*}
{where $\gamma^{-i\star}=(\gamma^{1\star},\ldots, \gamma^{i-1 \star},\gamma^{i+1\star},\ldots,\gamma^{N \star})$.}
Using Bayes' theorem, we can equivalently (for optimization purposes) rewrite the above optimization by
\begin{align}\label{eq:inf-ins}
    \inf_{u^i\in \mathbb{U}} \left\{g^i(H_1,u^{i})+g^{i}(H_2, u^{i})l^{i}\frac{P(H_2)}{P(H_1)}\right\},
\end{align}
where 
\begin{align}\label{eq:gi}
g^i(H_j,u^i):=\mathbb{E}^{\gamma^{-i\star}}\left[c(H,\gamma^{0}(u^{1:N}))\mid H=H_j\right]
\end{align}
for $j=1$ and $2$.
Since the expression inside the infimum in \eqref{eq:inf-ins} is linear in the likelihood ratio $l^i$, we conclude that an optimal response of the deviating sensor among all randomized policies belongs to the deterministic threshold type that only utilizes the likelihood ratio $l^i$.

Next, we state the optimality of the log-likelihood ratio test for the fusion center.

\begin{lemma}\label{lemma:2-PI}\cite{tsitsiklis1989decentralizedDet}
Consider $\mathcal{P}^{N}$ under Assumption \ref{Assump:density}. For any deterministic encoding policies $\gamma^{1:N}$ for each sensor that utilize $l^i$ instead of $y^i$, the optimal decoding policy of the fusion center is given by
\begin{align}\label{eq:h-test}
    \gamma^{0\star}(u^{1:N})=\begin{cases}
        1 & \text{if}\quad \Delta_{N}(u^{1:N}; \gamma^{1:N}) \geq t\\
        2 & \text{if}\quad \Delta_{N}(u^{1:N}; \gamma^{1:N}) < t
    \end{cases}
\end{align}
for some $t\in \mathbb{R}$, where $\Delta_{N}(u^{1:N}; \gamma^{1:N})$ is given by
\begin{align}\label{eq:DeltaN}
\Delta_{N}(u^{1:N}; \gamma^{1:N}):=\frac{1}{N} \sum_{i=1}^{N}\log \frac{\int P^{\gamma^i}(U^i=u^i|l^{i})P(dl^i|H_1)}{\int P^{\gamma^i}(U^i=u^i|l^{i})P(dl^i|H_2)}.
\end{align}
\end{lemma}

{We now explain the reasoning behind the optimality of the decoding policy presented in the above lemma. Following Lemma \ref{lem:TS-optimal}, $l^i$ are sufficient statistics for each sensor $i$. Focusing on deterministic encoding policies $\gamma^{1:N}$ for each sensor that utilize $l^i$ instead of $y^i$. Consider the fusion center's belief (posterior) on the hypothesis:
\begin{align*}
    P(H_i|u^{1:N})=\frac{P^{\gamma^{1:N}}(u^{1:N}\mid H_i) P(H_i)}{P(u^{1:N})}.
\end{align*}
The fusion center can make a decision based on the log ratio of the belief on $H_1$ and $H_2$ given by
\begin{align*}
    \log\left(\frac{P(H_1|u^{1:N})}{P(H_2|u^{1:N})}\right)=\log\left(\frac{P^{\gamma^{1:N}}(u^{1:N}\mid H_1)}{P^{\gamma^{1:N}}(u^{1:N}\mid H_2)}\right) + \log \frac{P(H_1)}{P(H_2)}
\end{align*}
which can be expressed as
\begin{align*}
 \log\left(\frac{P^{\gamma^{1:N}}(u^{1:N}\mid H_1)}{P^{\gamma^{1:N}}(u^{1:N}\mid H_2)}\right)&= \log\left(\prod_{i=1}^{N}\frac{P^{\gamma^{i}}(u^{i}\mid H_1)}{P^{\gamma^{i}}(u^{i}\mid H_2)}\right)\\
 &=  \sum_{i=1}^{N}\log\left(\frac{P^{\gamma^{i}}(u^{i}\mid H_1)}{P^{\gamma^{i}}(u^{i}\mid H_2)}\right)
\end{align*}
since conditioned on the hypothesis, sensors' actions are independent. Under Assumption \ref{Assump:density}, by considering the joint distribution on observations, we can show that the above is the same as $\Delta_{N}(u^{1:N}; \gamma^{1:N})$ without the normalization factor $1/N$. As a result, for any encoding policy of the sensors, the fusion center's decoding policy only depends on the above log-likelihood ratio term. }

Given the supporting results above, we present an equivalent problem that will be critical for our analysis. {For our threshold policies, $l^i$ and $\Delta_N$ (given in \eqref{eq:li} and \eqref{eq:DeltaN}) are sufficient statistics for each sensor and the fusion center, respectively. The expected cost $J(\gamma^{0:N})$ can be written as 
\begin{align}\label{eq:ctil}
    \int \tilde{c}_0(H,\Delta_N)\prod_{i=1}^{N} \mathbb{I}_{\{\gamma^i(l^i)\in du^i\}} P(dl^i|H)P(dH),
\end{align}
where $\tilde{c}_0(H,\Delta_N):=c\left(H,\gamma^0(\Delta_N(u^{1:N};\gamma^{1:N}))\right)$ is the new cost function.
Hence, following Lemma \ref{lem:TS-optimal} and Lemma \ref{lemma:2-PI}, we can reformulate the problem such that each sensor only has access to $l^i$ and the fusion center has access only to the likelihood ratio of $\Delta_N$ as depicted in Fig.~\ref{figure:2}.

\begin{figure}[ht]
    \centering
     \scalebox{.78}{
    \begin{tikzpicture}[node distance=5cm, auto]
        
        \node (E) [draw, rectangle, minimum width=1cm, minimum height=1cm, yshift=0.75cm] {$H$};
        \node (A) [draw, circle, minimum width=1cm, minimum height=2cm, right of=E, xshift=-2cm, yshift=0cm,text width=1.8cm] {Observation\\ \:\:\:channel};
        
        \node (B1) [draw, rectangle, minimum width=2cm, minimum height=1cm, right of=A, yshift=2cm] {Sensor $1$};
        \node (B2) [draw, rectangle, minimum width=2cm, minimum height=1cm, right of=A, yshift=0.5cm] {Sensor $2$};
        \node (Bdots) [draw=none, right of=A, yshift=-0.5cm] {$\vdots$};
        \node (BN) [draw, rectangle, minimum width=2cm, minimum height=1cm, right of=A, yshift=-1.5cm] {Sensor $N$};
        \node (C) [draw, circle, minimum width=2cm, minimum height=1cm, right of=Bdots, yshift=0.5cm, text width=1.8cm] {Observation\\ \:\:\:Channel};
        \node (D) [draw, rectangle, minimum width=2cm, minimum height=1cm, right of=C, yshift=0cm,xshift=-1cm] {Fusion Center};
        \node[text width=1cm] at (7,2.6) {$l^1$};
        \node[text width=1cm] at (7,1.5) {$l^2$};
        \node[text width=1cm] at (7,0) {$l^N$};
        \node[text width=1cm] at (9.8,2.6) {$u^1$};
        \node[text width=1cm] at (9.8,1.5) {$u^2$};
        \node[text width=1cm] at (9.8,0) {$u^N$};
        \node[text width=1cm] at (15.3,1.2) {$\Delta_N$};
        
        \draw[->] (E) -- (A);
        \draw[->] (A) -- (B1);
        \draw[->] (A) -- (B2);
        \draw[->] (A) -- (BN);
        \draw[->] (B1) -- (C);
        \draw[->] (B2) -- (C);
        \draw[->] (BN) -- (C);
        \draw[->] (C) -- (D);
    \end{tikzpicture}
    }
    \caption{The information structure under threshold encoding and log-likelihood ratio decoding policies. The middle column corresponds to the information of each sensor $i$, which is the likelihood ratio $l^i$, and the fusion center information is $\Delta_N$ computed from the sensors' actions.}\label{figure:2}
\end{figure}

Next, we show that an optimal solution exists in which the sensors' policies are of the threshold type. To this end, we require a topology for the space of threshold encoding policies. {We require a strong notion of convergence since the expected cost \eqref{eq:ctil} under the MAP decoding rule $\gamma^0$ of the fusion center is not continuous in the sensor's policies under the weak convergence topology.  If the encoding policies are viewed as joint probability measures on $\mathbb{L}\times \mathbb{U}$ for each sensor and are equipped with the weak convergence topology, then $\gamma^i_n$ converges weakly to $\gamma^i$ as $n\to \infty$ for every $i=1, \ldots, N$ does not, in general, imply continuity of the induced expected cost. In particular, the following convergence may fail:$$\inf_{\gamma^0\in \Gamma^0}\mathbb{E}^{\gamma^{1:N}_n}[ c(H,\gamma^0(u^{1:N})]\xrightarrow[]{n\to \infty}\inf_{\gamma^0\in \Gamma^0}\mathbb{E}^{\gamma^{1:N}}[ c(H,\gamma^0(u^{1:N})].$$ 
The reason is that the MAP rule $\gamma^0$ depends on both the joint message $u^{1:N}$ and the sensors’ encoding policies $\gamma^{1:N}$. Consequently, the expected cost is not the integral of a fixed continuous function; instead, the integrand itself varies with the encoding policies. A counterexample in the context of quantizers is given in \cite[Section 3.1.1]{yuksel2012optimization}, where it is shown that the expected cost need not be continuous under weak convergence of the observation channel. Therefore, weak convergence of encoding policies is not necessarily sufficient for establishing continuity of the expected cost, which is needed to prove the existence of an optimal solution. We thus require a stronger topology for policy convergence.}

In the following, we endow a topology for threshold policies using that introduced in  \cite{yuksel2012optimization} for quantizers utilizing the total variation distance. To this end, we view each threshold policy $\gamma^i$ as a stochastic kernel $\gamma^{i}(k|l^i)= \mathbb{I}_{\{l^i\in B_k^i\}}$
     for $k\in \mathbb{U}:=\{1, \ldots, M\}$, where $B_{k}^i=\left\{l^i\in \mathbb{L}| t_{k-1} < l^i \leq t_k\right\}$ are code cells. A sufficient statistic random variable $l^i$ is given by \eqref{eq:li}. 
In other words, a threshold policy $\gamma^i\in \mathcal{P}(\mathbb{U}|\mathbb{L})$ is a stochastic kernel, where $\mathbb{L}\subseteq \mathbb{R}$ is the set of all possible $l^i$. Equip $\mathcal{P}(\mathbb{U}|\mathbb{L})$ with the topology induced by the total variation distance $\|\cdot\|_{\sf TV}$ at any input measure $P$ on $\mathbb{L}$, that is, $\gamma^i_n$ converges to $\gamma^i$ as $n$ goes to infinity in total variation at input measure $P$ if as $n\to \infty$
\begin{align*}
    &\|P\gamma^i_n - P\gamma^i\|_{\sf TV}=\sup_{f:\|f\|_{\infty}\leq 1} \left| \int f(u,l)\gamma^i_n(du|l)P(dl) -\int f(u,l)\gamma^i(du|l)P(dl)\right|\to 0.
\end{align*}

We again denote the above set of threshold policies by $\Gamma_{\sf TS}^N$. {Since $\mathbb{U}$ is finite, the above topology is equivalent to
\[\sup_{\mathbb{A}\subseteq \mathbb{U}} \left|\int\gamma^i_n(\mathbb{A}|l)P(dl) - \int\gamma^i(\mathbb{A}|l)P(dl)\right|  \xrightarrow[n\to \infty]{}0\]
meaning that two policies are close at the input measure $P$ if they assign almost the same probability mass to every measurable event, uniformly.
We show that this makes the expected cost continuous in sensor policies, i.e., if $\gamma^i_n$ converges to $\gamma^i$ as $n\to \infty$ for all sensors in total variation, then
$$\inf_{\gamma^0\in \Gamma^0}\mathbb{E}^{\gamma^{1:N}_n}[ c(H,\gamma^0(u^{1:N})]\xrightarrow[n\to \infty]{}\inf_{\gamma^0\in \Gamma^0}\mathbb{E}^{\gamma^{1:N}}[ c(H,\gamma^0(u^{1:N})].$$}

Using \cite[Theorem 5.8]{yuksel2012optimization}, the following lemma shows that the space of threshold policies is compact under total variation. 
The proof is essentially from \cite{yuksel2012optimization}, but we include it in Appendix~\ref{app:lem-comptTV} for completeness.
\begin{lemma}\label{lem:compactness of TS}
    Under Assumption \ref{Assump:density}, $\Gamma_{\sf TS}^N$ is compact under total variation.
\end{lemma}

 The proof of the following theorem uses Lemma \ref{lem:compactness of TS} and is included in Appendix \ref{APP:the3}.

\begin{theorem}\label{the:thre}
Consider $\mathcal{P}^{N}$.  Let Assumption \ref{Assump:density} hold. Then, there exists an optimal solution $(\gamma^{0\star}, \gamma^{1\star}, \ldots, \gamma^{N\star})$; furthermore, this policy is a MAP rule for the fusion center and threshold policy for each sensor, i.e., $\gamma^{i\star}\in \Gamma_{\sf TS}^N$ for all $i\in \{1,\ldots, N\}$. 
\end{theorem}

The proof of the above theorem follows from the compactness of the space of threshold policies established in Lemma  \ref{lem:compactness of TS} and the continuity of the fusion center's optimization problem in threshold policies for the sensors. Compared to \cite{tsitsiklis1989decentralizedDet}, which discusses only the optimality of threshold policies, we establish the existence of an optimal encoding policy that is of threshold type. We note that, without such a threshold-type characterization of encoding policies (or quantizers, as in \cite{yuksel2012optimization}), the set of deterministic encoding policies is generally not weakly closed and hence not compact. Restricting to threshold-type policies also allows us to conclude total variation convergence under the Young topology, thereby establishing continuity of the expected cost. Together with compactness, this continuity guarantees the existence of an optimal policy.

Since one of our goals is to establish the existence and structural results for an optimal solution to limiting problems with an infinite number of sensors, we need to allow randomization among sensors' policies, which involves various relaxations given the decentralized information structure.

\section{Randomized Policies for Sensors and Information Structures for the Fusion Center}\label{sec:Randpolicies}.
Following the previous results, we focus on the equivalent formulation where the sensors only use $l^i$ as their new observation and the fusion center applies log-likelihood threshold policies that utilize only $\Delta_N$. From this and Assumption \ref{Assump:density}, we equip the space of randomized policies for each sensor separately with the weak-convergence topology. In particular, we focus on each {sensor}~$i$ separately and identify its randomized policy space $\Gamma^{i}$ via the set of probability measures on $\mathbb{U} \times \mathbb{L}$ with the fixed marginal $P$ on $\mathbb{L}$, i.e., 
\begin{align}\label{R-policy}
 &\Gamma^{i}:=\left\{\gamma^i \in \mathcal{P}(\mathbb{U}\times \mathbb{L}) \mid \gamma^i(\mathbb{U}, \cdot)=P(\cdot)\right\}
 \end{align}

To illustrate the endowed topology of $\Gamma^i$ and make it a Borel space, we define the set 
\begin{align}\label{thetaITanim}
 &\Theta^{i}:=\bigg\{\theta \in \mathcal{P}(\mathbb{U}\times \mathbb{L}) \bigg| \theta(\mathbb{A})=\int_{\mathbb{A}}\mathbb{I}_{\{g^{i}(l^{i})\in du^{i}\}}P(dl^{i}), g^i:\mathbb{L}\to \mathbb{U}\quad \forall \mathbb{A}\in \mathcal{B}(\mathbb{U}\times \mathbb{L})\bigg\},
 \end{align}
 where $\mathbb{I}$ denotes the indicator function. We equip the set $\Theta^{i}$ with the weak-convergence topology, in that the sequence $\{\theta^{i}_{n}\}_{n\geq 1}\subseteq \Theta^{i}$ converges to $\theta^{i}$ if and only if $\mathbb{I}_{\{g^{i}_{n}(l^{i})\in du^{i}\}}P(dl^{i})$ converges to $\mathbb{I}_{\{g^{i}(l^{i})\in du^{i}\}}P(dl^{i})$ weakly as $n\to \infty$. By \cite{BorkarRealization} (see also \cite[Chapter 4]{YukselBasar24}), $\Theta^{i}$ is the set of extreme points of $\Gamma^i$. This means that $\Gamma^i$ inherits Borel measurability and topological properties of the Borel measurable set $\Theta^{i}$. Hence, the sequence of randomized policies $\{\gamma^i_{n}\}_{n\geq 0}$ converges to a policy $\gamma^i$ if and only if $\gamma^i_{n}(du^i,dl^i):=\pi_{n}^{i}(du^i|l^{i})P(dl^{i})$ converges to $\gamma^i(du^i,dl^i):=\pi^{i}(du^i|l^{i})P(dl^{i})$ weakly as $n\to \infty$.

Since our analysis requires working with a compact and convex set of policies for the sensors, we need to allow a correlation in randomization among sensors. In the following subsection, we define various sets of randomized policies for the sensors utilized in our analysis and the main theorems.

\subsection{Finite Population Randomized Team Policies for Sensors}

Now that we have a standard Borel space formulation for $\Gamma^{i}$ for each sensor $i$, we can define the following set of randomized policies for the team of $N$ sensors as a set of Borel probability measures  $$L^{N}:=\mathcal{P}(\Gamma_{N})$$ with $\Gamma_{N}:=\Gamma^1 \times \ldots \times \Gamma^N$, where Borel $\sigma$-field $\mathcal{B}(\Gamma^{i})$ is induced by the topology defined before. We expand the fusion center's information to include the randomness utilized in the sensors' policies. 

Next, we define the set of randomized policies with common and independent randomness given by 
\begin{flalign*}
&\LCON:=\bigg\{P_{\pi} \in L^{N}\bigg{|}\forall \mathbb{A}_{i} \in \mathcal{B}(\Gamma^{i}): \nonumber\\
&P_{\pi}(\gamma^{1} \in \mathbb{A}_{1},\dots,\gamma^{N}\in \mathbb{A}_{N})=\int_{z\in [0,1]}\prod_{i=1}^{N}\widehat{P}_{\pi}^{i}(\gamma^{i}\in \mathbb{A}_{i}|z)\eta(dz), ~~~~\eta \in \mathcal{P}([0, 1])\bigg\},
\end{flalign*} 
where $\eta$ is the distribution of common but independent (from intrinsic exogenous system variables) randomness. Note that policies $\gamma^i$ are independent conditioned on $z$. In \cite{sanjari2020optimality}, it has been shown that $\LCON$ and $L^{N}$ are the same. Hence, the randomized policies in $L^{N}$ correspond to randomized policies induced by individual and common randomness. 

As mentioned in Section \ref{sec:Binary}, we need to allow the fusion center to have access to the realized sensors' policies unless there is a possible loss in performance; see our example provided in Section \ref{se:ExitExN}. Thus, we allow the fusion center to have access to the ({known common randomization}) information set $I^{\sf KCR}_N$ in \eqref{eq:IS-KR}
under which the common randomization device is known to the fusion center. Denote the set of all policies for the fusion center with this information set by $\Gamma^0_{\sf KCR}$.

For our analysis, we work with another information set for the fusion center which is amenable to our convergence argument as $N\to \infty$. For this information set, we allow the fusion center to observe, instead of randomization devices, the empirical measure of random policies $\gamma^{1:N}$ denoted by $\mu^N$, i.e.,
\begin{align*}
    \mu^N(\cdot)=\frac{1}{N}\sum_{i=1}^{N}\delta_{\gamma^i}(\cdot).
\end{align*}
Given $\mu^N\in \mathcal{P}(\Gamma^1)$, the fusion center can only determine the realized policies $\gamma^{1:N}$ of the sensors up to their permutation. In other words, the fusion center cannot distinguish the realized permutation of $\gamma^{1:N}$ by having only access to $\mu^N$. Lack of this knowledge might impact the performance, and hence, we additionally allow the fusion center to have access to the realized permutation by observing the realization of indices $I_{1}, \ldots, I_N$ (with $I_{i}\in\{1,\ldots,N\}$ valued) that allows the fusion center to infer the realized $\gamma^{1:N}$ for any randomized policy $P_{\pi}\in \LCON$. The above discussion yields that we need the information set of the fusion center to be 
\begin{align}\label{eq:IFS}
    I^{\sf FS}_N=\{u^{1:N}, \mu^N, I_{1}, \ldots, I_N\}.
\end{align}
Denote the set of all policies for the fusion center with this information set by $\Gamma^0_{\sf FS}$. Under the above information structure, the fusion center can infer realized $\gamma^{1:N}$. We show that given the above information structure for the fusion center, independent and common randomness does not improve the optimal expected cost, and hence, the above relaxation of the problem to sets of randomized policies $\LCON$ (which is identical to $L^N$) is a legitimate relaxation for team problems with $N$-Sensors. The proof of the following lemma is presented in Appendix \ref{app:leminfo}. 

\begin{lemma}\label{Lem:info-optimal-N}
Consider $\mathcal{P}^{N}$. Suppose that Assumption \ref{Assump:density} holds. Then,
\begin{align*}
    \inf\limits_{\gamma^{0}\in \Gamma^{0},\gamma^{1:N}\in \Gamma_N}J^N(\gamma^{0:N})= \inf\limits_{\gamma^0\in \Gamma^{0}_{\sf FS},P_{\pi}^{N} \in \LCON}J^N(\gamma^0, P_{\pi}^N).
\end{align*}
\end{lemma}

Before introducing a subset of $L^N$ including exchangeable policies, we formally introduce the concept of \emph{exchangeability} for finite and infinite collections of random variables and then policies.
\begin{definition}
Random variables $x^{1},x^{2},\dots,x^{N}$ defined on a common probability space are $N$-\it{exchangeable} if for any permutation $\sigma$ of the set $\{1,\dots,N\}$, 
\begin{align}
&{Pr}\bigg(x^{\sigma(1)} \in \mathbb{A}^{1},x^{\sigma(2)} \in \mathbb{A}^{2},\dots,x^{\sigma(N)}\in \mathbb{A}^{N}\bigg)={Pr}\bigg(x^{1} \in \mathbb{A}^{1},x^{2} \in \mathbb{A}^{2},\dots,x^{N}\in \mathbb{A}^{N}\bigg) \nonumber
\end{align}
for any measurable $\{\mathbb{A}^{1},\dots, \mathbb{A}^N\}$. Random variables $(x^{1},x^{2},\dots)$ are {\it infinitely-exchangeable} if they are $N$-\it{exchangeable} for all $N\in \mathbb{N}$.
\end{definition}

Using the above exchangeability notion for random variables, we introduce the set of exchangeable encoding policies for sensors given by
\begin{align}
\LEXN:=\bigg\{&P_{\pi} \in L^{N}\bigg{|}\forall \mathbb{A}_{i} \in \mathcal{B}(\Gamma^{i})~\text{and}~\forall\sigma \in S_{N}:\\
&P_{\pi}(\gamma^{1} \in \mathbb{A}_{1},\dots,\gamma^{N}\in \mathbb{A}_{N})=P_{\pi}(\gamma^{\sigma(1)} \in \mathbb{A}_{1},\dots,\gamma^{\sigma(N)}\in \mathbb{A}_{N})\bigg\}, \nonumber
\end{align}
where $S_{N}$ is the set of permutations of $\{1,\dots,N\}$. We note that $\LEXN$ is a subset of $L^{N}$, hence $\LCON$ because $L^{N}=\LCON$.

Since our interest lies in the existence of symmetric (identical) policies among sensors, we also introduce $\LCOSN$ as the set of symmetric policies with independent and common randomness:
\begin{eqnarray}
&&\LCOSN:=\bigg\{ P_{\pi} \in L^{N}\bigg{|}\forall\mathbb{A}_{i} \in \mathcal{B}(\Gamma^{i}): \nonumber\\
&& P_{\pi}(\gamma^{1} \in \mathbb{A}_{1},\dots,\gamma^{N}\in \mathbb{A}_{N})=\int_{z\in [0,1]}\prod_{i=1}^{N}\widehat{P}_{\pi}(\gamma^{i} \in \mathbb{A}_{i}|z)\eta(dz), \quad \eta \in \mathcal{P}([0, 1])\bigg\}, \nonumber
\end{eqnarray} 
where we drop the index $i$ on $\widehat{P}_{\pi}$ to capture the property that the independent randomization is identical through sensors. Also, define the set of randomized policies with only independent randomness as:
\begin{eqnarray}
L_{\sf PR}^{N}&:=&\left\{P_{\pi} \in L^{N}\mid \forall \mathbb{A}_{i} \in \mathcal{B}(\Gamma^{i}): P_{\pi}(\gamma^{1} \in \mathbb{A}_{1},\dots, \gamma^{N} \in \mathbb{A}_{N})  =\prod_{i=1}^{N}\widehat{P}_{\pi}^{i}(\gamma^{i}\in \mathbb{A}_{i})\right\}\nonumber.
\end{eqnarray} 
Finally, define the set of symmetric policies with independent randomness as follows
\begin{eqnarray}
\LPRSN&:=&\left\{P_{\pi} \in L^{N}\mid \forall \mathbb{A}_{i} \in \mathcal{B}(\Gamma^{i}): P_{\pi}(\gamma^{1} \in \mathbb{A}_{1},\dots, \gamma^{N} \in \mathbb{A}_{N})=\prod_{i=1}^{N}\widehat{P}_{\pi}(\gamma^{i}\in \mathbb{A}_{i})\right\}\nonumber.
\end{eqnarray} 

We note that $L_{\sf PR}^{N}\cap \LEXN=\LPRSN$, meaning that independent policies are exchangeable if they are symmetric.

Denote the space of independently randomized policies on deterministic threshold policies by $L_{\sf TS}^N\subset \LPR^N$. The following figure depicts the set relations between the space of policies for the finite population setting.

\begin{figure}[ht]
    \centering
    \scalebox{.8}{
    \begin{tikzpicture}
        \draw[thick] (-4,-3) rectangle (7,3);
        \node at (-2.8,2.5) {$L^N=\LCON$};
        
        \draw[thick] (0,0) circle(2.5);
        \draw[thick] (2,0) circle(2.5);
        \draw[ultra thick, color=blue!60] (2,0) circle(1.75);

        \node at (-1.2,1.5) {$\LEXN$};
        \node at (3.1,0.5) {\textcolor{blue}{$L_{\sf TS}^N$}};
        \node at (2,2.15) {$\LPR^N$};

        \begin{scope}
            \clip (0,0) circle(2.5);
            \shade[ball color=red!30,opacity=0.5] (2,0) circle(2.5);
        \end{scope}
        
        \node at (1,0) {\textcolor{purple}{$\LPRSN$}} ;
    \end{tikzpicture}
    }
    \caption{Sets of encoding policies for sensors. The relationship between sets of encoding policies for the sensors is illustrated in the figure: $L_{\sf PR}^{N}\cap \LEXN=\LPRSN$, $L_{\sf TS}^N\subset \LPR^N\subset \LCON$.}\label{Fig-policies}
\end{figure}

\subsection{Randomized Team Policies for Infinitely Many Sensors}

For the infinite population setting, we define sets $L, \LEX, \LCO, \LCOS, \LPR$, and $\LPRS$ using the Ionescu-Tulcea extension theorem by iteratively adding new coordinates to the probability measure. We define the set of probability measures $L$ on the infinite product Borel spaces $\Gamma=\prod_{i\in \mathbb{N}}\Gamma^{i}$ as 
\begin{eqnarray}
L:=\mathcal{P}(\Gamma)\label{eq:Linf}.
\end{eqnarray} 
Next we define randomized policies with common and independent randomization,
\begin{eqnarray}
&&\LCO:=\bigg\{P_{\pi} \in L\bigg{|}\forall\mathbb{A}_{i} \in \mathcal{B}(\Gamma^{i}): \nonumber\\
&&P_{\pi}(\gamma^{1} \in \mathbb{A}_{1},\gamma^{2} \in \mathbb{A}_{2},\dots)=\int_{z\in [0,1]}\prod_{i\in \mathbb{N}}\widehat{P}_{\pi}^{i}(\gamma^{i}\in \mathbb{A}_{i}|z)\eta(dz), \quad \eta \in \mathcal{P}([0, 1])\bigg\}. \nonumber
\end{eqnarray} 
 Note that $L_{\text{CO}}$ is a convex subset of $L$, and its extreme points are in the set of randomized policies with  independent randomness:
\begin{eqnarray}
\LPR&:=\bigg\{P_{\pi} \in L\bigg{|}\forall\mathbb{A}_{i} \in \mathcal{B}(\Gamma^{i}): P_{\pi}(\gamma^{1} \in \mathbb{A}_{1},\gamma^{2} \in \mathbb{A}_{2},\dots)=\prod_{i\in \mathbb{N}}\widehat{P}_{\pi}^{i}(\gamma^{i}\in \mathbb{A}_{i})\bigg\}\nonumber.
\end{eqnarray} 
For the above randomized policies, we expand the information set of the fusion center to
\begin{align}\label{eq:IFS-inf}
    I^{\sf KCR}=\{u^{1:\infty}, \gamma^{1:\infty}\}.
\end{align}
Next, we define the set of infinitely exchangeable policies as
\begin{eqnarray}
\LEX&:=&\bigg\{P_{\pi} \in L\bigg{|}\forall\mathbb{A}_{i} \in \mathcal{B}(\Gamma^{i}),~ \forall N\in\mathbb{N},~\text{and}~ \forall \sigma \in S_{N}: \nonumber\\
&&P_{\pi}(\gamma^{1} \in \mathbb{A}_{1},\dots,\gamma^{N}\in \mathbb{A}_{N})=P_{\pi}(\gamma^{\sigma(1)} \in \mathbb{A}_{1},\dots,\gamma^{\sigma(N)}\in \mathbb{A}_{N})\bigg\}. \nonumber
\end{eqnarray} 
Also, we symmetrically randomize policies with common and independent randomization as follows:
\begin{eqnarray}
&&\LCOS:=\bigg\{ P_{\pi} \in L\bigg{|}\forall \mathbb{A}_{i} \in \mathcal{B}(\Gamma^{i}): \nonumber\\
&& P_{\pi}(\gamma^{1} \in \mathbb{A}_{1},\gamma^{2}\in \mathbb{A}_{2}, \dots)=\int_{z\in [0,1]}\prod_{i\in \mathbb{N}}\widehat{P}_{\pi}(\gamma^{i} \in \mathbb{A}_{i}|z)\eta(dz), \quad \eta \in \mathcal{P}([0, 1])\bigg\}, \nonumber
\end{eqnarray} 
and further we define
\begin{eqnarray}
\LPRS&:=\bigg\{P_{\pi} \in L\bigg{|}\forall \mathbb{A}_{i} \in \mathcal{B}(\Gamma^{i}): P_{\pi}(\gamma^{1} \in \mathbb{A}_{1},\gamma^{2}\in \mathbb{A}_{2}, \dots)=\prod_{i\in \mathbb{N}}\widehat{P}_{\pi}(\gamma^{i} \in \mathbb{A}_{i})\bigg\}. \nonumber
\end{eqnarray} 

We now present a connection between $\LEX$ and $\LCOS$ using the de Finetti theorem (see, e.g.,\cite{aldous1981weak, Kallenberg}), which will establish that infinitely-exchangeable randomized policies are mixtures of conditionally i.i.d. randomized policies. We also recall a key result by Aldous \cite[Proposition 7.20]{aldous2006ecole}. We use these theorems for our main results as the number of sensors goes to infinity. These results were first introduced in  \cite{sanjari2020optimality} as a mathematical tool for establishing optimality of a symmetric and independent policy for mean-field team problems.

\begin{theorem}\label{the:defin}[A de Finetti Theorem for sensors individual policies]\cite{sanjari2020optimality}
Any infinitely-exchangeable randomized policy $P_{\pi}\in \LEX$ belongs to the set of randomized policies $\LCOS$ ($P_{\pi}\in \LCOS$), i.e., for any $P_{\pi}\in \LEX$, there exists a $[0, 1]$-valued random variable $z$ such that for any $\mathbb{A}_{i} \in \mathcal{B}(\Gamma^{i})$
\begin{flalign}
&P_{\pi}(\gamma^{1} \in \mathbb{A}_{1},\gamma^{2} \in \mathbb{A}_{2},\dots)=\int_{z\in [0, 1]}\prod_{i\in \mathbb{N}}\widehat{P}_{\pi}(\gamma^{i}\in \mathbb{A}_{i}|z)\eta(dz), ~~~~\eta \in \mathcal{P}([0, 1]),\label{eq:definiti}
\end{flalign} 
where $\widehat{P}_{\pi}(\cdot | z) \in {\cal P}(\Gamma^{i})$.
\end{theorem}

We note that $\LCOSN \subseteq \LEXN$; however, in general $\LCOSN \not= \LEXN$ (see a counter example in \cite{schervish2012theory}). 

\begin{theorem}\cite[Proposition 7.20]{aldous2006ecole}\label{the:ald}
Let $X:=(X_{1},X_{2},\dots)$ be an infinitely exchangeable sequence of random variables taking values in a Polish space $\mathbb{X}$. Suppose that $X$ is directed by a random measure $\alpha:\Omega \to \mathcal{P}(\mathbb{X})$, i.e.,
\begin{align}\label{mixture-iid}
    P(X_{1}\in A^1, X_{2}\in A^2, \ldots)= \int \prod_{i=1}^{\infty}\alpha(A^i) \eta(d\alpha) 
\end{align}
In other words, conditioned on $\alpha$, $(X_{1},X_{2},\dots)$ are i.i.d. random variables. 
Suppose that 
\begin{enumerate}
\item either for each $n$ $X^{(n)}=(X_{1}^{(n)},X_{2}^{(n)},\dots)$  is infinitely exchangeable directed by $\alpha_{n}$,
\item or $X^{(n)}=(X_{1}^{(n)},\dots, X_{n}^{(n)})$  is $n$-exchangeable with empirical measure $\alpha_{n}$.
\end{enumerate}
Then, $X^{(n)}$ converges in distribution to $X$ $\big(X^{(n)}\xrightarrow[n \to \infty]{{\sf d}}X\big)$  if and only if $\alpha_{n}\xrightarrow[n \to \infty]{{\sf d}}\alpha$ (by convergence in distribution to an infinitely exchangeable sequence, we mean the following: $X^{(n)} \xrightarrow[n \to \infty]{{\sf d}} X$ if and only if $(X_{1}^{(n)},\dots, X_{m}^{(n)})\xrightarrow[n \to \infty]{{\sf d}}(X_{1},\dots, X_{m})$ for each $m\geq 1$ \cite[page 55]{aldous2006ecole}).
\end{theorem}

In the following, we establish results on compactness and convexity of encoding policies, which will be used for our main theorems.  First, we have the following result on compactness and convexity of $\LEXN$ and $\LEX$ from \cite{sanjari2020optimality} (under the weak convergence topology).

\begin{lemma}[\cite{sanjari2020optimality}]\label{the:00}
Under Assumption \ref{Assump:density}, $\LEXN$ and $\LEX$ are convex and compact.
\end{lemma}

\section{Optimality of Exchangeable Policies for Finitely Many Sensor Setting}\label{se:ExitExN}

In the following, we present additional structural and existence results for optimal encoding policies. 
Recall that Lemma \ref{lemma:2-PI} guarantees that the fusion center policy is a MAP rule, and hence, without any loss of optimality, we can restrict our search for an optimal decoding policy to this class of policies under the information \eqref{eq:IFS}. For any randomized encoding policy of the sensors $P_{\pi}\in L^N=\LCON$, $\gamma^0$ will be of the form \eqref{eq:h-test}, where $\Delta_{N}(u^{1:N}; \gamma^{1:N})$ in  \eqref{eq:h-test} is replaced by $\overline{\Delta}_{N}(u^{1:N},P_{\pi},\mu^N,I_{1}, \ldots, I_N)$ given by
\begin{align}\label{eq:BaysDelta}
\overline{\Delta}_{N}(u^{1:N},P_{\pi},\mu^N,I_{1}, \ldots, I_N)=\frac{1}{N}\sum_{i=1}^{N}\log \frac{\int {\nu}(U^i=u^i|l^{i},P_{\pi},\mu^N,I_{1},\ldots, I_N)P(dl^i|H_1)}{\int {\nu}(U^i=u^i|l^{i},P_{\pi},\mu^N,I_{1},\ldots, I_N)P(dl^i|H_2)},
\end{align}
where 
\begin{align}
    {\nu}(U^i=u^i|l^{i},P_{\pi},\mu^N,I_{1},\ldots, I_N)=\gamma^{I_i}(U^i=u^i|l^i)
\end{align}
for the realized policy $\gamma^{I_i}$ of the sensor $i$.
The following lemma now establishes the exchangeability of an optimal encoding policy. The proof of the following lemma is included in Appendix \ref{APP:lem2}.

\begin{lemma}\label{lemma:1-EX}
Consider $\mathcal{P}^{N}$. Suppose that Assumption \ref{Assump:density} holds. Suppose that the fusion center's policy $\gamma^{0}$ is MAP given in \eqref{eq:h-test} with \eqref{eq:BaysDelta}. Then,  
\begin{itemize}
    \item [(i)] For any encoding policy $P_{\pi}\in L^N$,
    \begin{align}\label{eq:sigma-cost}
    J^N(\gamma^0, P_{\pi}^{\sigma})=J^N(\gamma^0, P_{\pi}),
\end{align}
where $\sigma$ is any permutation of $\{1,\ldots, N\}$ and $P_{\pi}^{\sigma}$ is a permutation of $P_{\pi}$, i.e., 
$P_{\pi}^{\sigma}(d\gamma^{1:N})=  P_{\pi}(d\gamma^{\sigma(1):\sigma(N)})$.
    \item [(ii)] Furthermore,
    \begin{align*}
    \inf_{P_{\pi}\in L^N} J^N(\gamma^0, P_{\pi})=\inf_{P_{\pi}\in \LEXN} J^N(\gamma^0, P_{\pi}).
\end{align*}
\end{itemize}
\end{lemma}

We now explain the reasoning behind the proof. Since observations $y^{1:N}$ are i.i.d. conditioned on the hypothesis, the ratios $l^i$ are also conditionally i.i.d. Using this together with the fact that $\mu^N$ is permutation invariant, we show that for any $P_{\pi}\in L^N$
\begin{align}\label{eq:exchang-Delta}
 \overline{\Delta}_{N}(u^{1:N},P_{\pi},\mu^N,I_{1}, \ldots, I_N)= \overline{\Delta}_{N}(u^{\sigma(1):\sigma(N)},P_{\pi}^{\sigma},\mu^N,I_{\sigma(1)}, \ldots, I_{\sigma(N)})
\end{align}
for all $u^{1:N}\in \mathbb{U}^{N}$. The expected cost can be written as
\begin{align*}
    J^N(\gamma^0, P_{\pi}^{\sigma})&=P\left(\cup \{\overline{\Delta}_{N}(u^{1:N},P_{\pi}^{\sigma},\mu^N,I_{\sigma(1)}, \ldots, I_{\sigma(N)})<t\}|H_1\right)P(H_1)\\
    &+P\left(\cup \{\overline{\Delta}_{N}(u^{1:N},P_{\pi}^{\sigma},\mu^N,I_{\sigma(1)}, \ldots, I_{\sigma(N)})\geq t\}|H_2\right)P(H_2).
\end{align*}
Consider the first term above. Following \eqref{eq:exchang-Delta}, $u^{1:N}$ leads to  $\overline{\Delta}_{N}(u^{1:N},P_{\pi},\mu^N,I_{1}, \ldots, I_N)<t$ if and only if $u^{\sigma(1):\sigma(N)}$ leads to $\overline{\Delta}_{N}(u^{\sigma(1):\sigma(N)},P_{\pi}^{\sigma},\mu^N,I_{\sigma(1)}, \ldots, I_{\sigma(N)})<t$. We also have $$P(\overline{\Delta}_{N}(u^{1:N},P_{\pi},\mu^N,I_{1}, \ldots, I_N)<t|H_1)=P(\overline{\Delta}_{N}(u^{\sigma(1):\sigma(N)},P_{\pi}^{\sigma},\mu^N,I_{\sigma(1)}, \ldots, I_{\sigma(N)})<t|H_1).$$ 
Using the convexity of the set $L^N$, then for any arbitrary policy ${P}_{\pi}$, we construct an exchangeable policy as an average of all its permutations. Since the realized policies can be deduced from $\mu^N$ and $I_1, \ldots, I_N$, an exchangeable policy attains the same performance as the average of the performance of the permutations of $P_{\pi}^{\sigma}$. Since by part (i), all the permutations attain the same performance, we conclude that for any given policy in $L^N$, we have an exchangeable policy in $\LEXN$ that attains the same performance.

Lemma \ref{lemma:2-PI} together with Lemma \ref{lemma:1-EX} yields that, under Assumption \ref{Assump:density}, without loss of optimality, the search for an optimal policy can be restricted to those that are $N$-exchangeable among sensors with a ratio-test-type decoding policy for the fusion center (which is permutation invariant). This, however, does not imply the existence of an optimal solution with this property.  Next, we use Theorem \ref{the:thre}, Lemma \ref{lemma:1-EX}, and Lemma \ref{lemma:2-PI} to show that an optimal solution exists and is exchangeable among sensors. The proof of the following theorem is included in Appendix \ref{APP:the2}.

\begin{theorem}\label{the:exi-N}
Consider $\mathcal{P}^{N}$ under Assumption \ref{Assump:density}. Then, there exists an optimal solution $(\gamma^{0\star},P_{\pi}^{\star})$, that is exchangeable among sensors with $\gamma^{0\star}$ of the form \eqref{eq:h-test} with \eqref{eq:BaysDelta}, i.e., $P_{\pi}^{\star}\in \LEXN$. 
\end{theorem}

We note that although Theorem \ref{the:thre} establishes that an optimal policy exists and is of threshold type for the sensors, and Theorem \ref{the:exi-N} establishes that an optimal policy exists and is exchangeable among the sensors, we cannot guarantee that such an optimal policy is simultaneously symmetric and of threshold type independently; see \cite{papastavrou1990decentralized} for an upper bound on the performance loss due to using identical threshold encoding policies (see e.g., \cite{DDglobal1987,reibman1987performance,sarma1983decentralized,tsitsiklis1988decentralized} for discussion on identical sensor policies). In the following, we show that we cannot strengthen the exchangeability condition to a symmetry condition without loss (in the finitely many sensor regime).

We expand the example provided in \cite{tsitsiklis1988decentralized}, showing that symmetric and independent policies need not be optimal for a finite number of sensors. {
\begin{example}\label{example1}
Suppose that $N=2$, $\mathbb{Y}=\{1,2,3\}$, and $\mathbb{U}=\{1,2\}$. Also, let $P(H_1)=P(H_2)=0.5$. Consider the following conditional distribution for observations:
\begin{align*}
    &P(y^i=1|H_1)=\frac{4}{5}, \qquad P(y^i=2|H_1)=\frac{1}{5}, \qquad P(y^i=3|H_1)=0,\\
    &P(y^i=1|H_2)=\frac{1}{3}, \qquad P(y^i=2|H_2)=\frac{1}{3}, \qquad P(y^i=3|H_2)=\frac{1}{3},
\end{align*}
for $i\in \{1,2\}$. For each sensor $i$, we compute the possible values $l^i$ using \eqref{eq:li} which are $\frac{5}{12},\frac{5}{3}$ and $\infty$ when $y^i=1,2$, and $3$, respectively.
Up to output label permutations and excluding constant quantizers, there are two nontrivial likelihood-ratio threshold policies:

\begin{itemize}
    \item Policy A: Let $\alpha$ be a real number in $(\frac{5}{12},\frac{5}{3})$. Then, $u^i=1$ if $l^i\in [0,\alpha]$  (that is $y^i=1$) and $u^i=2$ for $l^i\in [\alpha,\infty]$ (that is $y^i=2$ or $3$);
    \item Policy B: Let $\beta$ be a real number in $(\frac{5}{3},\infty)$. $u^i=1$ if  $l^i\in[0,\beta]$ (that is $y^i=1$ or $2$) and $u^i=2$ for $l^i\in[\beta,\infty]$ (that is $y^i=3$).
\end{itemize}

Suppose that sensors apply asymmetric policies first (for example, sensor 1 applies policy $A$ while sensor 2 applies policy $B$). Then, the optimal response of the fusion center obtained using \eqref{eq:DeltaN} is given by
\begin{align}\label{eq:FCP}
 u^{0}=\bigg\{\begin{array}{cc}
   1   &  \Delta_2(u^1,u^2)>0\\
    2  & \Delta_2(u^1,u^2)\leq 0
 \end{array}   
\end{align}
 We note that $\Delta_2(u^1,u^2)>0$ if and only if  $u^1=u^2=1$; otherwise $\Delta_2(u^1,u^2)\leq 0$. The optimal cost can be computed as follows:
\begin{align*}
    J(\gamma^{0\star:2\star})&=\frac{1}{2}\left\{P(\Delta_2(u^1,u^2)\leq 0|H_1)+P(\Delta_2(u^1,u^2)>0|H_2)\right\}\\
    &=\frac{1}{2}\left\{1-P(u^1=u^2=1|H_1)+P(u^1=u^2=1|H_2)\right\}\\
    &=\frac{19}{90}=0.21.
\end{align*}
Now, suppose that both sensors apply A. Then, an analogous analysis yields the cost $J(\gamma^{0\star:2\star})=\frac{53}{225}=0.24$. On the other hand, if they both apply policy B, then we get $J(\gamma^{0\star:2\star})=\frac{2}{9}=0.22$. This shows that asymmetric threshold policies can outperform symmetric ones for the finite population setting. 

Note that the symmetric independent sensor policy that assigns probability $0.5$ to $\gamma^i=A$ and $\gamma^i=B$ for each sensor is also suboptimal, with the total expected cost $\frac{1}{4}(\frac{2}{9}+\frac{53}{225}+\frac{38}{90})= 0.22$. However, an exchangeable policy that assigns a probability of $0.5$ to  $(\gamma^1=A,\gamma^2=B)$ and $(\gamma^1=B,\gamma^2=A)$ is also optimal since the fusion center observes $\mu^2=[1\quad 1]$ (meaning that one sensor selects policy $A$ and the other policy $B$) and given the indices of the sensor  $(I_1=1, I_2=2)$ or $(I_1=2, I_2=1)$ can infer the realized policy is $A, B$ or $B,A$, and hence, the performance is again $0.21$. Hence, an exchangeable policy exists, but not a symmetric and independent.

We note that if we consider a Bayesian formulation, where the fusion center only observes $u^{1:N}$ and randomized policy $P_{\pi}$, then again under the exchangeable policy that assigns a probability of $0.5$ to  $(\gamma^1=A,\gamma^2=B)$ and $(\gamma^1=B,\gamma^2=A)$, the optimal response of the fusion center is $u^{0}=1$ if and only if $u^1=u^2=1$ since only
\begin{align*}
    \log \frac{0.5[P^{A,B}(u^1=u^2=1|H_1)+P^{B,A}(u^1=u^2=1|H_1)]}{0.5[P^{A,B}(u^1=u^2=1|H_2)+P^{B,A}(u^1=u^2=1|H_2)]}=\log\frac{18}{5}>0
\end{align*}
and for other values of $u^1$ and $u^2$, we get a negative value. Hence, the performance remains the same. Although the same performance is achieved under the Bayesian fusion center in this example, this does not hold true in general. 

For example, suppose that both sensors observe the true hypothesis. Sensor $1$ selects a policy $C$ under which $P(u^1=1|H_1)=1$ and $P(u^1=2|H_2)=1$ (i.e., honest reporting the hypothesis) while Sensor $2$ selects a policy $D$ under which $P(u^2=1|H_1)=P(u^2=2|H_1)=0.5$ and $P(u^2=1|H_2)=P(u^2=2|H_2)=0.5$ (i.e., fully random reporting). 

We observe that in 
the setting where the fusion center knows the realized encoding policies, 
$\Delta_2(1,1)=\Delta_2(1,2)=\infty$ and $\Delta_2(2,2)=\Delta_2(2,1)=-\infty$. Hence, under the fusion's policy \eqref{eq:FCP}, the probability of the error becomes $0$, implying that the optimal performance is $0$. The permutation of this policy remains optimal as it again attains $0$. However, in 
the setting where the fusion center knows only the joint policy distribution under the exchangeable policy that assigns a probability of $0.5$ to  $(\gamma^1=C,\gamma^2=D)$ and $(\gamma^1=D,\gamma^2=C)$, we get $\Delta_2(2,1)=\Delta_2(1,2)=0$, and hence, the probability of error becomes 
\begin{align*}
   \frac{1}{2}\left(\frac{1}{2}P^{D,C}(u^1=2, u^2=1|H_1)+\frac{1}{2}P^{C,D}(u^1=1, u^2=2|H_1)\right)= \frac{1}{4},
\end{align*}
which is clearly not optimal. 
Hence, averaging the permutations of an optimal asymmetric policy need not preserve optimality if the fusion center does not know the realized policy assignments.
\end{example}
}

\section{Optimality of a Symmetric Independent Policy for Infinite Population Problem}\label{sec:Inf}

In this section, we show that an optimal policy for the decentralized detection problem with infinitely many sensors is symmetric, independent, and of threshold type. As mentioned before, since the probability of error goes to zero exponentially as the number of sensors goes to infinity, we study the error exponent as our new cost for the infinite population problem.

For finite $N$, Theorem \ref{the:exi-N} still holds under $J^N_{\sf EE}$. To study the asymptotic behavior of the error exponent, we need the following assumption from \cite{tsitsiklis1989decentralizedDet} that allows us to provide upper and lower bounds for the error exponent.

\begin{assumption}\label{Assump:inf-g}
For $j=1,2$ and $i\geq 1$, $\mathbb{E}\left[\left(\log l^i\right)^2\mid H_j\right]<\infty$.
\end{assumption}

In the following, we establish the existence of a symmetric (possibly randomized) encoding policy for the sensors. The proof of the following theorem is included in Appendix \ref{APP:the4}.

\begin{theorem}\label{the:Exist-inf}
Consider $\mathcal{P}^{N}_{\sf E}$ and $\mathcal{P}^{\infty}_{\sf E}$.     Let Assumptions \ref{Assump:density} and \ref{Assump:inf-g} hold. Then,
    \begin{itemize}
        \item [(i)] There exists an optimal solution $(\gamma^{0\star},P_{\pi}^{\star})$, that is symmetric and independent (possibly randomized) among sensors and with $\gamma^{0\star}$ of the form \eqref{eq:h-test}, i.e., $P_{\pi}^{\star}\in \LPRS$. 
        \item [(ii)] A symmetric encoding optimal solution $P_{\pi}^{R}\in \mathcal{P}(\Gamma^R)$ of
\begin{align}\label{eq:rep-encod}
 \inf_{P_{\pi}^{R}\in \mathcal{P}(\Gamma^R)} \min_{s\in [0,1]} \int \log\left( \sum_{u\in \mathbb{U}} g(H_2, u; \gamma^{R})^{1-s} g(H_1, u; \gamma^{R})^{s} \right) P_{\pi}^{R}(d\gamma^R), 
\end{align}
with
\begin{align}\label{eq:gamma-g-R}
    g(H, u; \gamma^{R})&:=P^{\gamma^{R}}(U^{R}=u|H)=\int{\gamma^{R}}(U^R=u|l^{R})P(dl^R|H)
\end{align}
is an optimal encoding policy for $\mathcal{P}^{\infty}_{\sf E}$. 
\end{itemize}
\end{theorem}
We now summarize the main steps of our proof technique. We first use the lower bound for the expected cost $J_{\sf EE}^{N}$ for any encoding policy from \cite[Lemma 1]{tsitsiklis1988decentralized} and \cite{shannon1967lower}, that is, for any encoding policy $\gamma^{1:N}$ of the sensors, we get
\begin{align*}
   \inf_{\gamma^{0}} &J^N_{\sf EE}(\gamma^{0:N})\geq \min_{s\in [0,1]}\frac{1}{N}\sum_{i=1}^{N}\log\left( \sum_{u^i\in \mathbb{U}} g(H_2, u^{i}; \gamma^{i})^{1-s} g(H_1, u^{i}; \gamma^{i})^{s} \right) \\
   &\left\{- \frac{1}{N} \sqrt{2\sum_{i=1}^{N}\frac{d^2}{ds^2}\log\left( \sum_{u^i\in \mathbb{U}} g(H_2, u^{i}; \gamma^{i})^{1-s} g(H_1, u^{i}; \gamma^{i})^{s} \right)}\bigg{|}_{s^{\star}_N}+\frac{\log(\min\{P(H_1),P(H_2)\}/2)}{N}\right\}.
\end{align*}
The second term on the right-hand side will vanish as $N\to \infty$ since the term with the second derivative is uniformly bounded under Assumption \ref{Assump:inf-g} following \cite[Proposition 3]{tsitsiklis1988decentralized}. Using this, we show that the limsup of the optimal expected cost $J_{\sf EE}^{N}$ under exchangeable encoding policies is the same as the limsup of that under infinitely exchangeable policies.  We note that Theorem \ref{the:exi-N} only allows us to restrict encoding policies to finite exchangeable policies (which might not be symmetric and independent). Here, we use a probabilistic argument based on a finite de Finetti theorem and a convergence result of an exchangeable sequence of random variables by Aldous \cite{aldous2006ecole} (see Theorem \ref{the:ald}) to show that the restriction can be refined to infinitely exchangeable encoding policies. This allows us to use the de Finetti theorem for policies in Theorem \ref{the:defin} to show the optimality of symmetric and independent sensors encoding policies, which, in turn, establishes the optimality of the representative encoding problem \eqref{eq:rep-encod}.

We now remark on the result of Theorem \ref{the:Exist-inf} and our analysis. 

\begin{remark}
First, we observe that in \eqref{eq:rep-encod}, the fusion center does not know the realized $\gamma^R$ and the corresponding private randomization device.  This is because, as $N$ goes to infinity, the empirical measures of $\mu^N$ converge to the law of symmetric policies. Hence, in the infinite population problem, under an optimal encoding policy of sensors that is symmetric and independent, the fusion center cannot infer the realized policies $\gamma^{R}$ (only the distribution $P_{\pi}^{\sf R}$ is known). Theorem \ref{the:Exist-inf} shows that asymptotically as $N\to \infty$, knowledge regarding the randomization device of the sensors for the fusion center does not improve the optimal performance. In fact, our analysis shows that a lack of knowledge about the randomization of sensors' policies does not incur any loss in performance for the detection problem asymptotically as $N$ goes to infinity. That is, as $N$ goes to infinity, the asymptotic optimal performance under the information set \eqref{eq:IFS-inf} where the realized sensors' policies are known to the fusion center is the same as that under the information set $I^0=\{u^{1:\infty}\}$, where the fusion center does not know the realization of policies.
\end{remark}

As a proposition to Theorem \ref{the:Exist-inf}, we formally state that an optimal policy exists (which consists of a MAP rule for the fusion center and an independent and identical policy for the sensors) for the infinite population setting and show that it is asymptotically optimal for large problems with error exponent cost. The proof of the following proposition is included in Appendix \ref{APP:Prop1}.

\begin{proposition}\label{the:thre-inf}
  Consider $\mathcal{P}^{N}_{\sf E}$ and $\mathcal{P}^{\infty}_{\sf E}$.  Let Assumptions \ref{Assump:density} and \ref{Assump:inf-g} hold.  Then,
  \begin{itemize}
      \item [(i)] There exists an optimal solution $(\gamma^{0\star}, \gamma^{R\star}, \gamma^{R\star},\ldots)$ for $\mathcal{P}^{\infty}_{\sf E}$ that is a MAP rule for the fusion center and an independent (possibly randomized) and also identical policy for each sensor. 
      \item [(ii)] A symmetric optimal solution $(\gamma^{0\star}, \gamma^{R\star}, \gamma^{R\star},\ldots)$ for $\mathcal{P}^{\infty}_{\sf E}$ restricted to $N$-sensors is $\epsilon_N$-optimal for $\mathcal{P}^{N}_{\sf E}$ with $\epsilon_N\to 0$ as $N\to \infty$.
  \end{itemize}
\end{proposition}

In \cite[Proposition 5.2 and Corollary 5.1]{tsitsiklis1989decentralizedDet}, the connection between finite population detection problems and their limit has been established. Here, in addition to establishing this connection, we justify that an optimal solution exists for the infinite population problem that is symmetric and independent. This, then, in turn, justifies the optimality of the detection problem with the representative sensor.

{
Before applying our approach and analysis to communication networks beyond decentralized detection, we first discuss how the results of this paper may extend to decentralized detection problems with imperfect communication channels.
\begin{remark}[Extension to imperfect communication channels]
    Our results in Sections \ref{sec:exitsN}--\ref{sec:Inf} extend to noisy communication channels, provided that the noise corrupting the sensors' messages is independent and identically distributed across sensors. 

Suppose that the fusion center observes a noisy message sent by the sensors: $$z^i=\kappa(u^i,v^i)$$ for i.i.d. noise $v^i$ such that $z^i\sim \mathcal{T}(\cdot|u^i)$ for a transition kernel $\mathcal{T}$ representing a noisy channel with a continuous conditional density $\tau(\cdot|u^i)$. Then, our results in Sections \ref{sec:exitsN}--\ref{se:ExitExN} hold with $g^i(H_j,u^i)$ in \eqref{eq:gi} replaced by 
\begin{align*}
g^i(H_j,u^i):=\int 
c\!\left(H_j,\gamma^0(z^i,z^{-i})\right)
\,\mathcal{T}(dz^i\mid u^i)
\prod_{k\neq i}
\left[
P(dy^k|H_j)\,
P^{\gamma^{k\star}}(du^k\mid y^k)\,
\mathcal{T}(dz^k\mid u^k)
\right],
\end{align*}
and \eqref{eq:DeltaN} replaced by
\begin{align*}
\Delta_{N}(z^{1:N}; \gamma^{1:N}):=\frac{1}{N} \sum_{i=1}^{N}\log \frac{\int \tau(z^i|u^i)P^{\gamma^i}(du^i|l^{i})P(dl^i|H_1)}{\int \tau(z^i|u^i)P^{\gamma^i}(du^i|l^{i})P(dl^i|H_2)}.
\end{align*}
Accordingly, the analysis can be adapted to this setting.

Our results in Section \ref{sec:Inf} also extend to noisy communication channels, replacing the bound in \eqref{eq:bound-shannon} to 
\begin{align*}
     \inf_{\gamma^{0}} &J^N_{\sf EE}(\gamma^{0:N})\geq \min_{s\in [0,1]}\frac{1}{N}\sum_{i=1}^{N}\log\left( \int g(H_2, z^{i}; \gamma^{i})^{1-s} g(H_1, z^{i}; \gamma^{i})^{s} \lambda(dz^i)\right) \\
   &- \frac{1}{N} \sqrt{2\sum_{i=1}^{N}\frac{d^2}{ds^2}\log\left( \int g(H_2, z^{i}; \gamma^{i})^{1-s} g(H_1, z^{i}; \gamma^{i})^{s} \lambda(dz^i)\right)}\bigg{|}_{s^{\star}_N}+\frac{\log(\min\{P(H_1),P(H_2)\}/2)}{N}.
\end{align*}
with  $\lambda$ is a Lebesgue measure and  
\begin{align*}
    g(H_2, z^{i}; \gamma^{i})=\int\tau(z^i|u^i){\gamma^{i}}(du^i|l^{i})P(dl^i|H_2).
\end{align*}
\end{remark}
}

\section{Robustness and Approximation for Finite and Infinite Sensor Problems}

{Our existence result in the previous sections provides more than the attainment of an optimal decentralized detection policy satisfying a structural property; it also supplies the compactness and continuity structure needed to study robustness properties of the detection problem under slight model uncertainty for the finite sensor problem and to obtain finite approximations for the limiting representative problem.

\subsection{Continuity and Robustness to Mismatch for Observation Channels}

In this subsection, we study robustness properties of the detection problem under model uncertainty. This is particularly important in decentralized detection, where the true observation channel \(Q\) for sensors is not known exactly. Commonly, the sensors and the fusion center design their decision rules using an estimated or approximating channel \(Q_n\), obtained from data or prior modeling assumptions. Therefore, it is important to determine whether small modeling errors lead only to small losses in performance. In this direction, we establish two fundamental properties using our framework and analysis in previous sections. First, we establish a continuity property: If an approximating observation channel \(\widetilde Q\) is close to the true channel \(Q\) in total variation, then the
corresponding optimal value under \(\widetilde Q\) remains close to the optimal value under \(Q\). Second, we establish a robustness property: If \(\widetilde\gamma^{0:N}\) is optimal for the inaccurate model \(\widetilde Q\), then
applying this policy to the true model \(Q\) still yields a performance close to the optimal performance achievable under \(Q\). These results show that the decentralized detection architecture is robust with respect to misspecification of the observation channel and justify the use of approximate models in the design of sensor and fusion policies. Such problems are well studied in the context of stochastic control problems in \cite{kara2020robustness, KaraYuksel2021Chapter}.

The following theorem establishes continuity and robustness for the finite-sensor problem. The proof is included in Appendix~\ref{app:rob}.

\begin{theorem}\label{the:Robust}
Consider $\mathcal{P}^N$. Suppose that $\mathbb{U}$ is a finite set. Let $Q_n$ be an inaccurate joint observation channel assumed by the sensors and the fusion center, and let $Q$ be the true one. Suppose that under $Q_n$:
\begin{itemize}
    \item [(a)] $y^1, \ldots, y^N$ are i.i.d., conditioned on $H$ with conditional distribution $\overline{Q}_n(\cdot|H)$.
\item [(b)] The conditional density $f^n_{y^i|H}$ of $y^i$ given $H$ for every $i\in\{1,\ldots,N\}$ exists,
$f^n_{y^i|H}(y^i|H_1)>0$ for all $y^i\in\mathbb{Y}$, and the likelihood ratio
\begin{equation}
 l^i_n = L_n(y^i) := \frac{f^n_{y^i|H}(y^i|H_2)}{f^n_{y^i|H}(y^i|H_1)} \qquad a.s. 
\end{equation}
is continuously differentiable in $y^i$ with a gradient that does not vanish, i.e., $\nabla L_n(y^i)\neq 0$ for all
$y^i\in\mathbb{Y}$.
\end{itemize}
The following two statements hold if $Q_n$ converges in total variation to $Q$ at input probability measure $P$ on $H_1$ and $H_2$ as $n\to \infty$, i.e., 
\begin{align}\label{eq:TV-Q}
    &\|PQ_n - PQ\|_{\sf TV}=\sup_{f:\|f\|_{\infty}<1} \left| \int f(y,H)Q_n(dy|H)P(dH) -\int f(y,H)Q(dy|H)P(dH)\right|\to 0.
\end{align}
\begin{itemize}
    \item [(i)] We have the continuity property:
\begin{align}\label{eq:Continuity}
     \left|\inf_{\gamma^{0:N}} J^N(\gamma^{0:N};Q_n)- \inf_{\gamma^{0:N}} J^N(\gamma^{0:N};Q)\right|\leq \|PQ_n-PQ\|_{\sf TV}\to 0
\end{align}

\item [(ii)] Robustness: Additionally, if ${\gamma}^{0\star:N\star}_n$ is an optimal policy under an inaccurate model ${Q}_n$, then 
    \[\left| J^N({\gamma}^{0\star:N\star}_n;Q)- \inf_{\gamma^{0:N}} J^N(\gamma^{0:N};Q)\right|\leq 2\|PQ_n-PQ\|_{\sf TV}\to 0.\]
\end{itemize}
\end{theorem}

We now provide a simple example where Theorem \ref{the:Robust} applies.
\begin{example}
Under the true observation channel \(Q\), assume that, conditioned on \(H\),
the observations are i.i.d. distributed by $\overline Q(dy^i\mid H)$.
Let \(\theta\neq 0\) such that 
\begin{equation*}
  y^i=  \begin{cases}
        w^i, & H_1\\
        \theta+w^i, & H_2
    \end{cases}
\end{equation*}
with \(w^1,\ldots,w^N\) i.i.d. \(N(0,1)\). This means $\overline Q(\cdot\mid H_1)=N(0,1)$ and $
\overline Q(\cdot\mid H_2)=N(\theta,1)$. Now  define an
approximating independent channel \(Q_n\) that converge by $Q_n(dy^1,\ldots,dy^N\mid H)
=
\prod_{i=1}^N \overline Q_n(dy^i\mid H)$,
where 
\begin{equation*}
  y^i=  \begin{cases}
        w^i, & H_1\\
        \theta_n+w^i, & H_2
    \end{cases}
\end{equation*}
i.e., $\overline Q_n(\cdot\mid H_1)=N(0,1)$ and 
$\overline Q_n(\cdot\mid H_2)=N(\theta_n,1)$.
If \(\theta_n\to\theta\), we have pointwise convergence of the conditional
densities:
\[
f^n_{y^i\mid H}(y\mid H_j)
\to
f_{y^i\mid H}(y\mid H_j),
\qquad j=1,2.
\]
By Scheffe's theorem,
$\|\overline Q_n(\cdot\mid H_j)-\overline Q(\cdot\mid H_j)\|_{\sf {TV}}
\to 0$ for 
$j=1,2$. Since the channels are finite products, it follows that
\[
\|Q_n(\cdot\mid H_j)-Q(\cdot\mid H_j)\|_{\sf {TV}}
\to 0,
\qquad j=1,2,
\]
which implies \eqref{eq:TV-Q}.
This framework can be used for a data-driven approach: During a training phase, suppose that we collect labeled observations $Z_1,\ldots,Z_n \sim \overline Q(\cdot\mid H_2)=N(\theta,1)$,
independently. Define the data-driven estimator
\[
\theta_n:=\frac1n\sum_{k=1}^n Z_k .
\]
By the strong law of large numbers,
$\theta_n\to\theta$ a.s. as $n\to \infty$.

\end{example}

Now, we discuss the implications of our robustness result when the true observation channel $Q$ has conditionally correlated observations. Unfortunately, a conditionally correlated observation channel cannot, in
general, be approximated in total variation by conditionally independent
channels with vanishing error (that is because if $\{Q_n\}_n$ is a conditionally independent sequence of channels for all $n\geq 1$ and converges in the sense of \eqref{eq:TV-Q} to a channel $Q$, then $Q$ must also be conditionally independent). However, a conditionally correlated observation channel may still be approximated by such channels up to a small nonzero error.

\begin{example}
Let the true
observation channel be
$
Q(\cdot\mid H_j)=N(m_j,\Sigma_\rho)$,
for $j=1,2$,
where $m_1=0$, $
m_2=\theta \mathbf 1$,
and $\Sigma_\rho
=
(1-\rho)I_N+\rho \mathbf 1\mathbf 1^\top$ for 
$0<\rho<1$.
Thus the observations are conditionally correlated with $\operatorname{Cov}(Y^i,Y^k\mid H_j)=\rho$,
for $i\neq k$.

Now approximate \(Q\) by the independent channel
$\widetilde{Q}(\cdot\mid H_j)=N(m_j,I_N)$ for $j=1,2$. Note that we cannot find an approximate independent channels \(Q_n\) that converge to \(Q\) in the total variation unless
\(\rho\to0\), but for small \(\rho\), $\widetilde{Q}$ is close to \(Q\) in the total variation as we show below. By Pinsker's inequality
\begin{align*}
\|Q(\cdot\mid H_j)-\widetilde{Q}(\cdot\mid H_j)\|_{\mathrm{TV}}
&\leq
\sqrt{2
D_{\mathrm{KL}}\bigl(N(m_j,\Sigma_\rho)\,\|\,N(m_j,I_N)\bigr)}\,.
\end{align*}
Since \(\operatorname{tr}(\Sigma_\rho)=N\), we obtain
\begin{align*}
D_{\mathrm{KL}}\bigl(N(m_j,\Sigma_\rho)\,\|\,N(m_j,I_N)\bigr)&=\frac12
\left(
\operatorname{tr}(\Sigma_\rho)-N-\log |\Sigma_\rho|
\right)\\
&=
-\frac12\log |\Sigma_\rho|\\
&=-\frac{1}{2}\log\left((1-\rho)^{N-1}\bigl(1+(N-1)\rho\bigr)\right).
\end{align*}
Since
\[
\begin{aligned}
-\log|\Sigma_\rho|
&=
-(N-1)\log(1-\rho)
-\log\bigl(1+(N-1)\rho\bigr),
\end{aligned}
\]
using the Taylor expansions
$
\log(1-\rho)
=
-\rho-\frac{\rho^2}{2}+O(\rho^3)$,
and
$
\log\bigl(1+(N-1)\rho\bigr)
=
(N-1)\rho
-{(N-1)^2\rho^2}/{2}
+O(\rho^3)$,
we obtain
\begin{align*}
-\log|\Sigma_\rho|
&=
(N-1)\rho+\frac{N-1}{2}\rho^2
-(N-1)\rho
+\frac{(N-1)^2}{2}\rho^2
+O(\rho^3) \\
&=
\frac{N(N-1)}{2}\rho^2
+O(\rho^3).
\end{align*}
Therefore, $\|Q(\cdot\mid H_j)-\widetilde{Q}(\cdot\mid H_j)\|_{\sf{TV}}$ and hence, $\|PQ-P\widetilde{Q}\|_{\sf{TV}}$ for small \(\rho\) can be shown to be of order $O(\rho\sqrt{N(N-1)})$.
Thus, for fixed \(N\), the mismatch between the correlated and independent observation models vanishes as \(\rho \to 0\).
\end{example}

\subsection{Approximation for the Limiting Representative Problem}

In Proposition \ref{the:thre-inf}, we show that an optimal policy for the infinite sensor problem, obtained by solving \eqref{eq:rep-encod}, yields a near-optimal policy for the large finite-sensor problem with respect to the error-exponent cost. This policy has a particularly useful structure: The fusion center uses a MAP rule, while the sensors use an independent and identical representative policy. Hence, in the large-population regime, the
original decentralized detection problem can be reduced to the design of a single representative sensor rule. This structural reduction is significant because it removes the dependence of the policy search space on the number of sensors. Instead of optimizing over all decentralized policies for \(N\) sensors, one only needs to optimize over a single representative encoder. Nevertheless, the resulting representative
problem \eqref{eq:rep-encod} is still computationally nontrivial. The objective
depends nonlinearly on the induced message distributions and the optimization over the representative encoder is generally non-convex, as already noted in \cite{tsitsiklis1988decentralized}. Therefore, although the infinite-population formulation provides a scalable structural characterization, it does not by itself give a directly computable solution method. To address this issue, we next construct a finite approximation of the limiting representative problem \eqref{eq:rep-encod}. The main idea is to quantize the likelihood-ratio space and restrict attention to representative encoders that are constant on the resulting quantization cells. This converts the representative optimization problem into a finite optimization problem, which can in principle be solved by exhaustive search. We then show that, as the quantization becomes finer, the optimal value of the finite approximation converges to the optimal value of the original representative problem.

Consider the limiting representative problem in \eqref{eq:rep-encod} under Assumptions \ref{Assump:density} and  \ref{Assump:inf-g}. Suppose that $\mathbb{Y}$ is compact. Then, since the likelihood function is continuous by Assumption \ref{Assump:density}, $\mathbb{L}$ is compact and hence is totally bounded. This means that we can find a sequence $(\{l_{n,i}\}_{i=1}^{k_n})_{n\geq 1}$ of finite grids in $\mathbb{L}$ such that
\begin{align*}
   \sup_{l\in \mathbb{L}} \min_{i=1, \ldots, k_n} d_{\mathbb{L}}(l,l_{n,i})<\frac{1}{n} \qquad \forall n\geq 1, 
\end{align*}
where $k_n$ is the size of the finite grid and $d_{\mathbb{L}}$ is a metric on $\mathbb{L}$. We define the finite space $\widehat{\mathbb{L}}_n:=\{l_{n,1},\ldots,l_{n,k_n}\}$, which is the quantization of $\mathbb{L}$ under nearest-neighbor quantizer $q_n$ given by
\begin{align*}
q_n(l)=\argminA_{l_{n,i}\in \widehat{\mathbb{L}}_n}\:d_{\mathbb{L}}(l,l_{n,i})\quad \forall l\in \mathbb{L}.
\end{align*} 
Let  $A_{n,i}=\{l\in \mathbb{L}\:|\:q_n(l)=l_{n,i}\}$. Define for any $P\in \mathcal{P}(\mathbb{L}|\mathbb{H})$  $$\widehat{P}_{n}(\cdot|H)=\sum_{i=1}^{k_n}P(A_{n,i}\mid H)\delta_{l_{n,i}}(\cdot) \qquad \forall n.$$ Since $\mathbb{L}$ is compact, for any $P\in \mathcal{P}(\mathbb{L}|\mathbb{H})$, there exists a sequence $\{\widehat{P}_n\}_n$ with $\widehat{P}_n\in\mathcal{P}(\widehat{\mathbb{L}}_n|\mathbb{H})$  that converges weakly to $P$ as $n\to \infty$.

Now, suppose that the observation channel $\overline{Q}(dy^R|H_i)$ for $i=1,2$ is non-atomic. Then, under Assumption \ref{Assump:density}(iii), the conditional distribution of the likelihood values $P(dl^R|H_i)$ is also non-atomic for $i=1,2$. That is because
$P(\{a\}\mid H_j)
=
\overline Q\bigl(L^{-1}(\{a\})\mid H_j\bigr)$, and
under Assumption \ref{Assump:density}(iii), $\{y\in\mathbb Y|L(y)=a\}$
has zero \(\overline Q(\cdot| H_j)\)-measure, and hence, $P(\{a\}| H_j)=0$ for all $a\in\mathbb L$. Let 
\begin{align*}
    \mathcal{G}=\left\{(g(H_i, u; \gamma^{R}))_{i=1,2, u\in \mathbb{U}}\:\: \text{s.t.}\: g(H_i, u; \gamma^{R})=\int{\gamma^{R}}(U^R=u|l^{R})P(dl^R|H_i) \:\bigg|\:\gamma^R\in \Gamma^R\right\}.
\end{align*}
Let
$\Phi_n^R$ be the set of all deterministic encoders $\phi_n^R:\mathbb{L}\to\mathbb{U}$ that are constant on the quantization
cells $l^R\in A_{n,k}$. Then, define $\mathcal{G}^n$ by replacing the randomized encoder $\gamma^R$ with $\phi_n^R$, i.e.,
\begin{align*}
    \mathcal{G}^n=\left\{(g(H_i, u; \phi^{R}_n))_{i=1,2, u\in \mathbb{U}}\:\: \text{s.t.}\:g(H_i, u; \phi^{R}_n)=\int \mathbb{I}_{\{\phi^{R}_n(l^{R})=u\}}{P}(dl^R|H_i)\:\bigg|\:\phi^{R}_n\in \Phi_n^R \right\}.
\end{align*}
Then, since $\mathbb{U}$ is finite, following the denseness argument in \cite[Proof of Theorem 5.1]{yuksel2023borkar} under Young topology for encoding policies (see also \cite[Theorem 2.7]{beiglbock2018denseness}, \cite[Theorem 3]{milgrom1985distributional}, \cite{balder1997consequences}, \cite{borkar1988probabilistic}, or \cite[Chapter 7]{castaing2004young}), we get that for any $g\in \mathcal{G}$, we can find a sequence $\{g_n\}_n$ with  $g_n\in \mathcal{G}_n$ that converges to $g$ as $n\to \infty$. Let
\[
S_m:=\left\{0,\frac1m,\ldots,\frac{m-1}{m},1\right\}
\]
for some $m\geq 1$. Therefore, we define an approximate representative encoder problem 
\begin{align}\label{eq:rep-encod-approx}
J_{\sf EE}^{\star, n,m}
:=
\inf_{\phi_n^R\in\Phi_n^R}
\min_{s\in S_m}
\log\left(
\sum_{u\in\mathbb U}
g(H_2,u;\phi_n^R)^{1-s}
g(H_1,u;\phi_n^R)^s
\right).
\end{align}
In the following theorem, we show that as the quantization becomes finer, the optimal value $J_{\sf EE}^{\star, n,m}$ converges to the true one. The proof of the following theorem is included in Appendix \ref{APP:approx-proof}.

\begin{theorem}\label{the:Approx-rep}
    Consider the limiting representative problem in \eqref{eq:rep-encod} under Assumptions \ref{Assump:density} and \ref{Assump:inf-g}. Suppose additionally that $\mathbb{Y}$ is compact and the observation channels $\overline{Q}(dy^i|H_j)$ for $j=1,2$ are atomless for all $i\in \mathbb{N}$. Then, 
\begin{align*}
  \lim_{n\to \infty}  \lim_{m\to\infty}J^{\star,n,m}_{\sf EE}=J^{\star}_{\sf EE},
\end{align*}
where $J^{\star}_{\sf EE}$ corresponds to the optimal value in \eqref{eq:rep-encod}.
\end{theorem}
The above theorem implies that an optimal solution of the finite optimization problem in \eqref{eq:rep-encod-approx} is near-optimal for large $n$ and $m$. The problem in \eqref{eq:rep-encod-approx} can be solved by exhaustive search over $|\mathbb{U}|^{k_n}(m+1)$. 

We now give two simple numerical examples. 

\begin{enumerate}
    \item We illustrate the finite approximation procedure on a binary hypothesis testing
problem with scalar Gaussian observations. Under the two hypotheses, the
representative observation satisfies
\[
H_1:\quad Y^R\sim \mathcal N(0,1),
\qquad
H_2:\quad Y^R\sim \mathcal N(1,1).
\]
Let $\mathbb U=\{0,1,2\}$. Fig.~\ref{fig:numerical-convergence} demonstrates the approximation for the representative encoding problem and also the convergence of the approximation value under the quantized optimal policy of the representative encoder as the number of sensors grows. Note that in this example $\mathbb{Y}$ is not compact, but the approximation method using quantization can also be established for non-compact settings, following \cite{SaYuLi15c}.
\item For another simulation example when $\mathbb{Y}$ is compact, we take $\mathbb{Y}=[0,1]$ and let
\[ H_1:\quad Y^R\sim \mathrm{Unif}[0,1],\qquad H_2:\quad  f_{Y^R\mid H_2}(y) = \frac{\alpha e^{\alpha y}}{e^\alpha-1}, \qquad 0\le y\le 1, \] with \(\alpha=1.25\). The results are shown  
in Fig.~\ref{fig:numerical-convergence}, where we use $\mathbb U=\{0,\ldots,9\}$.
\end{enumerate}

\begin{figure}[ht]
    \centering
    \includegraphics[width=0.8\textwidth]{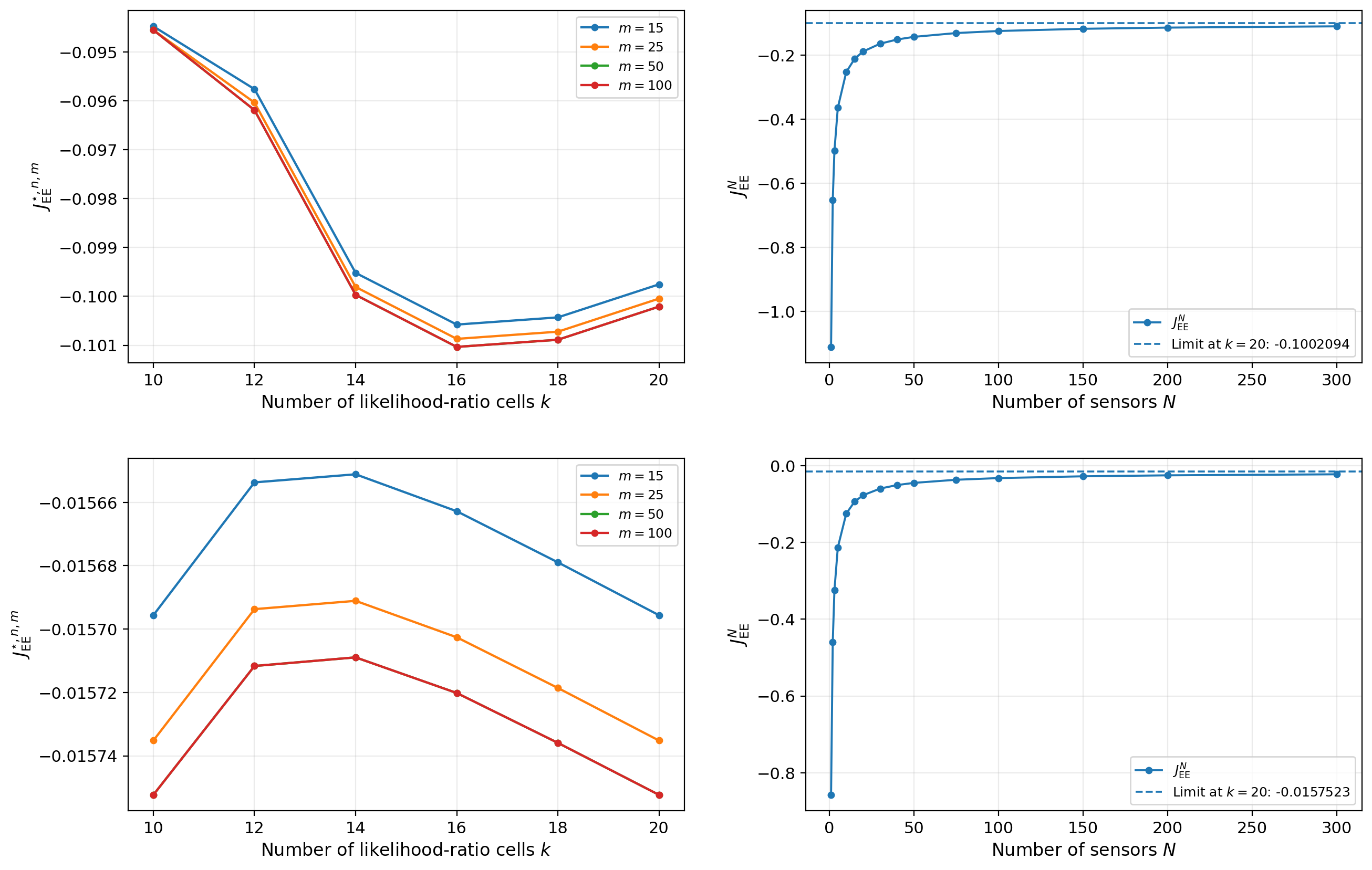}
    \caption{{Left: convergence of the finite approximation \(J_{\sf EE}^{\star,n,m}\) as the number of likelihood-ratio cells \(k_n\)
    and the \(s\)-grid size \(m\) increase.
    Right: finite-population error-exponent cost
    \(J^N_{\sf EE}\) versus $N$ 
    under the representative encoder obtained with \(k_n=20\) and \(m=100\).
    The dashed line is the representative limit for the fixed policy. The top panels correspond to the first example, while the bottom panels correspond to the second example.}}
    \label{fig:numerical-convergence}
\end{figure}

\section{Concluding Remarks and Future Work}\label{sec:Conc}

\subsection{Conclusion}

We studied a class of decentralized binary detection problems involving a single fusion center and many sensors. In the finite population setting, we established the existence and optimality of an independent threshold-based encoding policy for the sensors and a MAP rule for the fusion center. Additionally, we showed that an exchangeable encoding policy is also optimal. In the infinite population setting, we proved the existence and optimality of an independent, symmetric encoding policy under the error exponent cost. Using this, we demonstrated the optimality of a representative sensor equivalent problem and provided an approximation for the large but finite sensor regime.

{
\subsection{The exchangeability approach beyond decentralized detection: Multi-access channels}

The method developed in this paper has potential implications beyond decentralized detection. As a preliminary step in this direction, we briefly study a simple multi-access communication model motivated by massive machine-type communications (mMTC) and Internet of Things (IoT) networks, where a large number of devices (sensors) send data to a fusion center \cite{bockelmann2016mMTC}. In some studies on unsourced random access (URA), such as \cite{shao2020cooperative} and \cite{cakmak2023joint}, each device in a group of spatially distributed devices contributes partial and noisy data to the inference of a global hypothesis (state), as in the decentralized detection problem considered in this paper. Hence, the theoretical results on optimality and exchangeability also apply to such settings. Detailed investigation of the applications of the results in this paper to wireless communication networks employing URA \cite{ozates2024unsourced} is a possible direction for future work.

To show that the exchangeability and symmetry framework developed in the paper applies not only to decentralized detection but also to broader application areas, we consider a simple random-access communication model in this section under a heavy-traffic assumption. This illustrative extension is not intended as a comprehensive treatment of multi-access communication networks; rather, it highlights how permutation invariance and large-population limits can again lead to symmetric independent policies.

In this problem, $N$ users (or, in the limit as $N \to \infty$, infinitely many users) repeatedly decide whether to retransmit a packet over a broadcast channel under decentralized information. A popular algorithm, known as ALOHA, requires transmitters to use randomized policies based on local feedback. 

We consider an $N$-user (or transmitter) slotted random-access collision channel. At each time $t$, user $i$ observes its local queue length $L_t^i$ and chooses an action
\[
U_t^i \in \{0,1\},
\]
where $U_t^i=0$ denotes no transmission and $U_t^i=1$ denotes transmission. The queue length of user $i$ evolves according to
\begin{align}\label{eq:Li}
L_{t+1}^{i}
=
\max\bigg(0,
L_t^{i}+A_t^{i}
-
1_{\{C_t=1,\;U_t^{i}=1\}}1_{\{L_t^{i}\neq 0\}}
\bigg),
\qquad i=1,\dots,N,
\end{align}
where $A_t^i$ denotes the number of exogenous packet arrivals to user $i$ during slot $t$, and $1_{\mathbb{A}}$ denotes the indicator function of an event $\mathbb{A}$. The variable $C_t$ indicates whether the transmission attempt in slot $t$ is successful, and is given by
\[
C_t:=1_{\{\sum_{i=1}^N U_t^{i}=1\}}.
\]
Thus, a packet is successfully transmitted in slot $t$ if and only if exactly one user transmits in that slot. If no user transmits, the channel is idle; if two or more users transmit simultaneously, a collision occurs, and no packet is successfully delivered. Therefore, the departure term for user $i$ is equal to one only when user $i$ transmits, its queue is nonempty, and no other user transmits (see \eqref{eq:Li}).

We assume that the initial queue lengths $\{L_0^i\}_{i=1}^N$ are i.i.d. across users. We focus on the heavy-traffic regime in which the queues are asymptotically nonempty for large $t$.  Additionally, we assume that $A_t^{i}$ is an independent Bernoulli random variable with 
\begin{align}
    \mathbb{P}(A_t^{i}=1)= \frac{a}{N} \quad \forall i=1, \ldots, N.
\end{align}
for some $\bar\alpha<a$ and for all $t\geq 0$ with $ 1\leq\bar\alpha<\infty$.

User $i$ acts under decentralized information, having access to the history of $L_{t}^{i}$, transmitting with probability $\lambda^{N}_{i}$ given by
    \begin{align}\label{eq:pol-lambda}
    \lambda^{N}_{i}
    :=
    \frac{\alpha^{N}_i}{N}
    \end{align}
    for some $\alpha_i^N\in [0, \min\{N,\bar\alpha\}]$.
Under a (jointly) randomized policy ${\pi}^N$ of users as a joint distribution on $(\alpha^{N}_1, \ldots, \alpha^{N}_N)$, the objective is to arrive at a policy profile where the mismatch of transmission and data availability is minimized asymptotically as $t\to \infty$. For each $N$, the objective is to maximize the finite population performance criterion, which is expected throughput
    \[
    J^N_t(\pi^N)
    :=
    \mathbb E^{\pi^N}
    \left[R_t^N(L_{t}^{1},U_{t}^{1}, \ldots, L_{t}^{N},U_{t}^{N})
    \right]
    \]
for large $t$, where the per-slot reward is given by
\[
R_t^N(L_{t}^{1},U_{t}^{1}, \ldots, L_{t}^{N},U_{t}^{N})
:=
\sum_{i=1}^N 1_{\{U^i_t=1\}}
1_{\{L_t^{i}\neq 0\}}
\prod_{j\neq i}(1-1_{\{U^{j}_t=1\}}).
\]
In particular, when $t$ is large under the heavy-traffic assumption, we have $L_t^{i}> 0$ for all $i\in \mathbb{N}$, and hence, we instead focus on a static problem where the expected cost for a policy $\pi^N$ is given by
 \[
T^N(\pi^N):=\mathbb{E}^{\pi^{N}}\!\left[\sum_{i=1}^N 1_{\{U^i=1\}} \prod_{j\neq i} (1-  1_{\{U^j=1\}})\right].
\]

In the following theorem, we show that, for large $N$, a near-optimal policy is identical and independent (each user sends its packets according to an identical and independent policy). Our analysis parallels that used in the previous section for decentralized detection, thereby highlighting the usefulness of a stochastic-team formulation and the importance of probabilistic methods based on exchangeability in the study of communication-network problems. The proof of the following theorem is provided in Appendix~\ref{APP:Thm6}.

\begin{theorem}\label{the:Applic}
Consider the random-access problem above. Then:
\begin{itemize}
    \item [(i)] For finite $N$, there exists an optimal policy $\pi^{N\star}$ that is exchangeable among users.
    \item [(ii)] For large $N$, there exists an asymptotically optimal policy $\pi^{N\star}$ that is symmetric and independent across users. In particular, the asymptotic optimal policy $\pi^{N\star}$ has a transmission probability 
\[
\lambda_i^{N\star}=\frac{\alpha^i}{N}, \qquad i=1,\dots,N
\]
with $\alpha^i$ are independent with identical distribution $\delta_{1}(\cdot)$ and attains the expected throughput $e^{-1}$ approximately for large $t$ and $N$, i.e., 
\[\lim_{N\to \infty}\lim_{t\to \infty}J_t^N(\pi^{N\star})=e^{-1}.\]
\end{itemize}
\end{theorem}
We now outline our proof strategy. Existence and exchangeability of an optimal policy in Part (i) follows from an analogous argument used in the proofs of Lemma \ref{lemma:1-EX} and Theorem \ref{the:thre} by first showing the space of policies is compact under total variation, and the expected reward is continuous, then by constructing an exchangeable policy as an average of all its permutations, since for finite $N$ the throughput and dynamics are exchangeable, there exists an optimal policy ${\pi}^{N\star}$ which is exchangeable. For Part (ii), we first find an asymptotic upper bound for the optimal expected throughput. Then, we show that an infinitely exchangeable policy is optimal by attaining the upper bound. Then, we show that the asymptotic throughput of an infinitely exchangeable policy depends on its directing measure only through its mean and that the upper bound is attained when the common randomness is degenerate at $1$, leading to asymptotic optimality of the independent identical policy $\alpha^i\sim \delta_{1}(\cdot)$.

We make the following remark about the results above. The class of policies considered in \eqref{eq:pol-lambda} is restricted to low transmission probabilities for large \(N\), since \(\alpha_i^N\) takes values in a compact set. This restriction is important because, if arbitrary transmission probabilities were allowed, then for any \(N\) in heavy traffic, an expected throughput of $1$ could be achieved by any joint randomized policy that assigns probability one to configurations of  $\{\alpha_i^N\}_{i=1}^N$ in which one component equals \(N\) and all others equal \(0\). Equivalently, only one user transmits a packet with probability 1, while all other users do not. Such joint randomization effectively allows coordination among users. By contrast, the low-transmission-probability regime in \eqref{eq:pol-lambda} limits joint randomization and prevents sufficient coordination among users when \(N\) is large. In fact, as Theorem \ref{the:Applic} suggests, for large \(N\), coordination through joint randomization is not necessary in the sense of approximate optimality.

In our analysis, allowing joint randomization across policies would break down the exchangeability argument, since that argument requires the expected cost to depend only on the individual actions/policies and the empirical measure. Under low transmission probabilities, we show that the expected throughput can be upper-bounded solely in terms of the \(\alpha_i^N\)'s and their empirical measures, not their joint distribution.}

\appendix

\subsection{Proof of Lemma \ref{lem:TS-optimal}}\label{APP:lem TS}
 The optimality of threshold policies in the above theorem follows a similar idea to that in  \cite{tsitsiklis1989decentralizedDet}. Fix policies of sensors with index $k\in \{1, \ldots, i-1, i+1, \ldots, N\}$ to policies $\gamma^{k\star}$.  Consider the deviating sensor's optimization problem 
\begin{align}
 \inf_{u^i\in \mathbb{U}} \left\{\sum_{j=1}^{2} P(H_j|y^i)\mathbb{E}^{\gamma^{-i\star}}\left[c(H_j,\gamma^{0}(u^{1:N}))\mid H_j\right]\right\}.
\end{align}
By Bayes' theorem, almost surely
\begin{align*}
    \frac{P(H_2|y^i)}{P(H_1|y^i)}=L(y^{i}) \frac{P(H_2)}{P(H_1)}.
\end{align*}
Hence, the deviating sensor's optimization problem can be written as follows:
\begin{align*}
    \inf_{u^i\in \mathbb{U}} \left\{g^i(H_1,u^{i})+g^{i}(H_2, u^{i})l^{i}\frac{P(H_2)}{P(H_1)}\right\},
\end{align*}
where 
\begin{align*}
    g^i(H_j,u^i):=\mathbb{E}^{\gamma^{-i\star}}\left[c(H,\gamma^{0}(u^{1:N}))\mid H=H_j\right].
\end{align*}
Hence, the optimal response encoding policy for the sensor becomes 
\begin{align}
    \gamma^{i\star}(l^{i})\in\begin{cases}
        \argminA\limits_{u^i\in \mathbb{U}} \left\{g^i(H_2,u^{i})\right\} & \text{if}\: l^i=\infty\\
        \argminA\limits_{u^i\in \mathbb{U}} \left\{g^i(H_1,u^{i})+g^{i}(H_2, u^{i})l^{i}\frac{P(H_2)}{P(H_1)}\right\} & \text{if}\: l^i<\infty
    \end{cases}
\end{align}
 In the second case, for each $u^i\in\mathbb{U}$ we have an affine function of $l^i$ with the slope of
$g^i(H_2,u^i)\frac{P(H_2)}{P(H_1)}$. The minimum of finitely many affine functions is concave and piecewise affine,
and the set of values of $l^i$ for which a given $u^i$ attains the minimum is an interval.  Hence, we can conclude that a threshold policy $\widehat{\gamma}^i \in \Gamma_{\sf TS}^i$ for the sensor $i$ attains the minimum above since permutation of the action for the thresholds does not change the cost. Hence, the restriction to threshold policies for sensors is without loss of optimality.

\subsection{Proof of Lemma \ref{lem:compactness of TS}}\label{app:lem-comptTV}
The proof proceeds in two steps. In Step 1, we establish the compactness of the space of threshold policies for each sensor. In Step 2, we show that the product space of threshold policies is also compact in total variation.

\begin{itemize}

\item  [{\bf Step 1.}]{\bf Compactness of threshold policy space for each sensor under total variation.} In this step, following \cite[Lemma 5.2]{yuksel2012optimization}, we show that the space of threshold policies for each sensor is compact in total variation at any input measure $P$. 

Let $\nu= P \times \lambda$, where $\lambda$ is the counting measure on $\{1, \ldots, M\}$. For any measurable set $B\subset \mathbb{R}^n$ and $k=1, \ldots, M$, $\nu(B\times \{k\})=P(B)\lambda(\{k\})=P(B)$. Hence,
\begin{align*}
    P\gamma^i(\{k\} \times B)=P(B \cap B_{k}^{i})\leq P(B) = \nu(B\times \{k\})
\end{align*}
Since any measurable set $D\in \mathbb{L} \times \mathbb{U}$ can be written as a disjoint union of $D_k \times \{k\}$ with $D_{k}=\{l^i \in \mathbb{L}|(l^i,k)\in D\}$, we conclude that $P\gamma^i(D)\leq \nu(D)$. Following \cite[Lemma 4.1]{bogachev2007measure}, we conclude that threshold policies for each sensor are relatively compact in total variation at any input measure $P$. Therefore, for any sequence $\gamma^i_n$, there exists a subsequence that converges to $\gamma^i$. Next, we show that the limit policy $\gamma^i$ is also a threshold policy.

We use \cite{yuksel2012optimization} and \cite{gyorgy2003codecell} to show that if an input measure $P$ is absolutely continuous with respect to the Lebesgue measure, then the space of threshold policies is closed in total variation at any input measure $P$. 

The code cells $B_{k}^i=\left\{l^i| t_{k-1} < l^i \leq t_k\right\}$ are convex sets. This is because for any $l^i_1$ and $l^i_2$ elements of $B^i_k$, we can observe that $\alpha l^i_1 + (1 - \alpha) l^i_2\in B^i_k$ for $\alpha \in [0,1]$. Since the code cells are convex, if the input measure $P$ is absolutely continuous with respect to the Lebesgue measure, then \cite{gyorgy2003codecell} implies that a subsequence $\gamma^i_n$ converges to $\gamma^i$ in the sense that 
\begin{align}
    \lim_{n\to \infty} P(B^{i}_{n,k} \Delta B^i_k)=0
\end{align}
     for all $s\in \{1, \ldots M\}$, where $B^{i}_{n,s}$ and $B^{i}_{s}$ are code cells under $\gamma^i_n$ and $\gamma^i$, respectively. Hence, following the analysis in \cite{yuksel2012optimization}, we can show that

\begin{align*}
    \|P\gamma^i_n - P\gamma^i\|_{\sf TV}
    &=\sup_{f:\|f\|_{\infty}<1} \left| \sum_{s=1}^{M}\int f(s,l)\gamma^i_n(s|l)P(dl) -\int f(s,l)\gamma^i(s|l)P(dl)\right|\\
    &=\sup_{f:\|f\|_{\infty}<1} \left| \sum_{s=1}^{M}\int f(s,l)\left(\mathbb{I}_{l\in B^{i}_{n,s}} - \mathbb{I}_{l\in B^{i}_{s}}\right)P(dl)\right|\\
    &\leq \sup_{f:\|f\|_{\infty}<1}  \sum_{s=1}^{M}\int_{B^{i}_{n,s} \Delta B^i_s} |f(s,l)|P(dl)\\
    &\leq \sum_{s=1}^{M}P(B^{i}_{n,s} \Delta B^i_s) \to 0
\end{align*}
as $n\to \infty$. This justifies that the space of threshold policies is closed in total variation at any input measure $P$ since, following Assumption \ref{Assump:density}(iii), $P$ admits a density, and hence, it is absolutely continuous with respect to the Lebesgue measure. 

\item  [{\bf Step 2.}]{\bf Compactness of the joint threshold policy space.} In this step, we show that the joint space of policies for sensors is compact in total variation. Since the observations are independent, conditioned on the hypothesis $H$, we have
\begin{align*}
    &\|P\prod_{i=1}^{N}\gamma^i_n - P\prod_{i=1}^{N}\gamma^i\|_{\sf TV}\\
&=\sup_{f:\|f\|_{\infty}<1} \left| \int f(u^{1:N},l^{1:N})\prod_{i=1}^{N}\gamma^i_n(du^i|l^i)P(dl^i|H) -\int f(u^{1:N},l^{1:N})\prod_{i=1}^{N}\gamma^i(du^i|l^i)P(dl^i|H)\right|\\
&\leq\sup_{\widehat{f}:\|\widehat{f}\|_{\infty}<1} \left| \int \widehat{f}(u^{1},l^{1}, H)\gamma^1_n(du^1|l^1)P(dl^1|H) -\int \widehat{f}(u^{1},l^{1}, H)\gamma^1(du^1|l^1)P(dl^1|H)\right|\\
&+\sup_{\widehat{f}:\|\widehat{f}\|_{\infty}<1} \bigg| \int \widehat{f}(u^{2:N},l^{2:N},H)\prod_{i=2}^{N}\gamma^i_n(du^i|l^i)P(dl^i|H) -\int \widehat{f}(u^{2:N},l^{2:N},H)\prod_{i=2}^{N}\gamma^i(du^i|l^i)P(dl^i|H)\bigg|\\
&= \|P\gamma^1_n - P\gamma^1\|_{\sf TV}\\
&+\sup_{\widehat{f}:\|\widehat{f}\|_{\infty}<1} \bigg| \int \widehat{f}(u^{2:N},l^{2:N},H)\prod_{i=2}^{N}\gamma^i_n(du^i|l^i)P(dl^i|H) -\int \widehat{f}(u^{2:N},l^{2:N},H)\prod_{i=2}^{N}\gamma^i(du^i|l^i)P(dl^i|H)\bigg|,
\end{align*}
where
\begin{align*}
    \widehat{f}(u^{1},l^{1}, H)&=\int f(u^{1:N},l^{1:N})\prod_{i=2}^{N}\gamma^i_n(du^i|l^i)P(dl^i|H)\\
    \widehat{f}(u^{2:N},l^{2:N}, H)&=\int f(u^{1:N},l^{1:N})\gamma^1(du^1|l^1)P(dl^1|H).
\end{align*}
Repeating this procedure, we conclude that $\|P\prod_{i=1}^{N}\gamma^i_n - P\prod_{i=1}^{N}\gamma^i\|_{\sf TV}$ converges to zero as $n\to \infty$.
\end{itemize}

\subsection{Proof of Theorem \ref{the:thre}}\label{APP:the3}

Suppose each sensor applies threshold policies $\gamma^{1:N}$. Thanks to Lemma \ref{lem:compactness of TS}, we only need to show that the infimum of the expected cost is continuous in the sensors' policies. We have

\begin{align*}
   &\bigg| \inf_{\gamma^0\in \Gamma^0} \int c(H,\gamma^{0}(u^{1:N})) \prod_{i=1}^{N}\gamma^{i}_n(du^{i}|l^i)P(dl^i|H_k)P(dH_k) - \\& \inf_{\gamma^0\in \Gamma^0}\int c(H,\gamma^{0}(u^{1:N})) \prod_{i=1}^{N}\gamma^{i}(du^{i}|l^i)P(dl^i|H_k)P(dH_k)\bigg|\\
   &\leq \sup_{\gamma^0} \bigg|  \int c(H,\gamma^{0}(u^{1:N})) \prod_{i=1}^{N}\gamma^{i}_n(du^{i}|l^i)P(dl^i|H_k)P(dH_k) - \\&\int c(H,\gamma^{0}(u^{1:N})) \prod_{i=1}^{N}\gamma^{i}(du^{i}|l^i)P(dl^i|H_k)P(dH_k)\bigg|\\
   &\leq\left\|P\prod_{i=1}^{N}\gamma^i_n - P\prod_{i=1}^{N}\gamma^i\right\|_{\sf TV}\to 0
\end{align*}
as $n\to \infty$.

Following Lemma \ref{lemma:2-PI}, for any threshold policies $\gamma^{1:N}$, the optimal response of the fusion center $\gamma^0$ is a ratio test. Then, the optimality of threshold policies for sensors in Lemma \ref{lem:TS-optimal} completes the proof.

\subsection{Proof of Lemma \ref{Lem:info-optimal-N}} \label{app:leminfo} 
Let $\gamma^{0}$ be a MAP rule. We have 
   \begin{align}
    \inf\limits_{\gamma^{0}\in \Gamma^{0}_{\sf FS},\gamma^{1:N}\in \Gamma_N}J^N(\gamma^{0:N})&= \inf\limits_{\gamma^{1:N}\in \Gamma_{N}}J^N(\gamma^{0\star},\gamma^{1:N})\label{eq:lem4-1}\\
    &=\inf\limits_{P_{\pi}^{N} \in \LCON}J^N(\gamma^{0\star}, P_{\pi}^N)\label{eq:lem4-2}\\
    &=\inf\limits_{\gamma^0\in \Gamma^{0}_{\sf FS},P_{\pi}^{N} \in \LCON}J^N(\gamma^0, P_{\pi}^N)\label{eq:lem4-3}
\end{align}
where \eqref{eq:lem4-1} follows from Lemma \ref{lemma:2-PI}. Fixing a policy of the fusion center to the MAP rule $\gamma^{0\star}$, \eqref{eq:lem4-2} follows from \cite[Theorem 2.1(ii)]{YukselSaldiSICON17}) since $\LCON$ is a convex set with the set of extreme points in those policies that are deterministic conditioned on a common randomness (which is known to the fusion center under \eqref{eq:IFS}), and since the expected cost is linear in randomized encoding policies. The last equality \eqref{eq:lem4-3} follows from the fact that under \eqref{eq:IFS}, the fusion center can infer the common randomness and realized randomized policies $\gamma^{1:N}$, and hence, MAP rules remain optimal for the fusion center as the encoding policies can be without loss to be considered deterministic conditioned on a common randomness.
\subsection{Proof of Lemma \ref{lemma:1-EX}}\label{APP:lem2}
\begin{itemize}
\item [Part (i):] Since observations $y^{1:N}$ are i.i.d. conditioned on the hypothesis, the ratios $l^i$ are also conditionally i.i.d. Consider any $P_{\pi}\in L^N$. Using this together with the fact that $\mu^N$ is permutation invariant
\begin{align*}
&\int \prod_{i=1}^{N} \nu(u^{\sigma(i)}|l^{i},P_{\pi}^\sigma,\mu^N,I_{\sigma(1)},\ldots, I_{\sigma(N)})P(dl^i|H_1)\\
&=\int \prod_{i=1}^{N} \nu(u^{\sigma(i)}|l^{\sigma(i)},P_{\pi}^\sigma,\mu^N,I_{\sigma(1)},\ldots, I_{\sigma(N)})P(dl^{\sigma(i)}|H_1)\\
&=\int \prod_{i=1}^{N} \nu(u^i|l^{i},P_{\pi},\mu^N,I_{1},\ldots, I_N)P(dl^i|H_1)\\
&=\int \prod_{i=1}^{N}\gamma^{I_{i}}(u^i|l^i)P(dl^i|H_1).
\end{align*}
This implies that for any $P_{\pi}\in L^N$
\begin{align}\label{eq:pf-exchang-Delta}
 \overline{\Delta}_{N}(u^{1:N},P_{\pi},\mu^N,I_{1}, \ldots, I_N)= \overline{\Delta}_{N}(u^{\sigma(1):\sigma(N)},P_{\pi}^{\sigma},\mu^N,I_{\sigma(1)}, \ldots, I_{\sigma(N)})
\end{align}
for all $u^{1:N}\in \mathbb{U}^{N}$.

The expected cost can be written as
\begin{align*}
    J^N(\gamma^0, P_{\pi}^{\sigma})&=P\left(\cup \{\overline{\Delta}_{N}(u^{1:N},P_{\pi}^{\sigma},\mu^N,I_{\sigma(1)}, \ldots, I_{\sigma(N)})<t\}|H_1\right)P(H_1)\\
    &+P\left(\cup \{\overline{\Delta}_{N}(u^{1:N},P_{\pi}^{\sigma},\mu^N,I_{\sigma(1)}, \ldots, I_{\sigma(N)})\geq t\}|H_2\right)P(H_2).
\end{align*}
Consider the first term above. Following \eqref{eq:pf-exchang-Delta}, $u^{1:N}$ leads to  $\overline{\Delta}_{N}(u^{1:N},P_{\pi},\mu^N,I_{1}, \ldots, I_N)<t$ if and only if $u^{\sigma(1):\sigma(N)}$ leads to $\overline{\Delta}_{N}(u^{\sigma(1):\sigma(N)},P_{\pi}^{\sigma},\mu^N,I_{\sigma(1)}, \ldots, I_{\sigma(N)})<t$. We also have $$P(\overline{\Delta}_{N}(u^{1:N},P_{\pi},\mu^N,I_{1}, \ldots, I_N)<t|H_1)=P(\overline{\Delta}_{N}(u^{\sigma(1):\sigma(N)},P_{\pi}^{\sigma},\mu^N,I_{\sigma(1)}, \ldots, I_{\sigma(N)})<t|H_1)$$ since 
\begin{align*}
&\int \prod_{i=1}^{N} \nu(U^i=u^{\sigma(i)}|l^{i},P_{\pi}^\sigma,\mu^N,I_{\sigma(1)}, \ldots, I_{\sigma(N)})P(dl^i|H_1)\\
&=\int \prod_{i=1}^{N}\gamma^{I_{\sigma(i)}}(u^{\sigma(i)}|l^{\sigma(i)})P(dl^{\sigma(i)}|H_1)\\
&=\int \prod_{i=1}^{N}\gamma^{I_{i}}(u^i|l^i)P(dl^i|H_1)\\
&=\int \prod_{i=1}^{N} \nu(U^i=u^i|l^{i},P_{\pi},\mu^N,I_{1},\ldots, I_N)P(dl^i|H_1).
\end{align*}
An analogous argument for the case of conditioning on $H_2$ leads to \eqref{eq:sigma-cost}.

\item [Part (ii):] Using the convexity of the set $L^N$, then for any arbitrary policy ${P}_{\pi}$, we construct an exchangeable policy as an average of all its permutations, i.e., $$\overline{P}_{\pi}=\frac{1}{|S_N|}\sum_{\sigma \in S_N} {P}_{\pi}^{\sigma}.$$ Since the realized policies can be deduced from $\mu^N$ and $I_1, \ldots, I_N$, an exchangeable policy attains the same performance as the average of the permutations of $P_{\pi}^{\sigma}$, i.e.,
\begin{align*}
J^N(\gamma^0,\overline{P}_{\pi})=\frac{1}{|S_N|}\sum_{\sigma \in S_N} J^N(\gamma^0,{P}_{\pi}^{\sigma}).
\end{align*}
Since by part (i), all the permutations attain the same performance, we conclude that for any given policy in $L^N$, we have an exchangeable policy in $\LEXN$ that attains the same performance.
\end{itemize}

\subsection{Proof of Theorem \ref{the:exi-N}}\label{APP:the2}
By Theorem \ref{the:thre}, there exists an optimal detection policy belonging to the threshold type among all policies in $\LCON$, which is the same as $L^N$. This implies the existence of an optimal detection policy among $L^N$. By Lemma \ref{lemma:1-EX}, we can conclude that an optimal policy is without loss exchangeable, and this completes the proof. 
\subsection{Proof of Theorem \ref{the:Exist-inf}}\label{APP:the4}
The proof of part i proceeds in three steps. In Step 1, we show that a sequence of optimal exchangeable sensors' policies and their empirical measures admits a subsequence that converges to an infinitely exchangeable policy and its directing measure. In Step 2, we study the limit of the optimal exponent cost as the number of sensors goes to infinity and show that the limsup equals the representative cost.
In Step 3, we show that
the limiting policy of the optimal exchangeable $N$-sensor policy under a subsequence, which is symmetric and independent among sensors, is optimal for the exponent problem with an infinite number of sensors. 
\begin{itemize}
    \item [{\bf Step 1.}] {\bf Convergence property of optimal exchangeable policies and their empirical measure.}
Following Lemma \ref{the:00} and a finite de Finetti theorem in \cite[Theorem 13]{diaconis1980finite} (using an analogous argument as that in the proof of \cite[Lemma 2]{sanjari2020optimality}), by considering the indexes of optimal exchangeable optimal policies as a sequence of i.i.d. random variables with uniform distribution on the set $\{1,\dots,N\}$, we can construct a sequence of  infinitely exchangeable policies close in total variation (for any finite marginal) to an exchangeable optimal policies of the sensors, which converges to a limiting policy that is infinitely exchangeable.  This yields that a subsequence of exchangeable optimal policies of the sensors $\gamma^{1\star:n\star}_{n}$ converges in distribution (for any finite marginal) to an infinitely exchangeable policy $\gamma^{1\star:\infty\star}$. Following Theorem \ref{the:ald}, this implies that the subsequence of empirical measures $\mu^{n\star}$ of exchangeable optimal policies $\gamma^{1\star:n\star}_{n}$ of the sensors converges to the directing measure $\mu^{\star}$ of the limiting infinitely exchangeable policy $\gamma^{1\star:\infty\star}$. Hence, we conclude that a subsequence $P_{\pi}^{n\star}$ induced by exchangeable optimal policies converges to $P_{\pi}^{\star}\in \LEX$ and the empirical measure $\mu^{n\star}$ converges to the directing measure of $\mu^{\star}$ of $P_{\pi}^{\star}$.

Let $u^{1\star:n\star}_{n}$ be actions induced by exchangeable optimal policies $\gamma^{1\star:n\star}_{n}$. Since observations $y^{1:n}$ are i.i.d. under the change of measure argument, so are $l^{1:n}$. Hence, we can conclude that $u^{1\star:n\star}_{n}$ are exchangeable. Along the same reasoning, we get that $u^{1\star:n\star}_{n}$ converge in distribution to infinitely exchangeable actions $u^{1\star:\infty\star}$ induced by $\gamma^{1\star:\infty\star}$. By Theorem \ref{the:ald}, their empirical measures of  $u^{1\star:n\star}_{n}$ also converge weakly to the directing measure of the infinitely exchangeable sequence $u^{1\star:\infty\star}$. Let 
\begin{align}\label{eq:lambda-n}
\Lambda_n(\cdot)&:=\frac{1}{n} \sum_{i=1}^{n} \delta_{\{\gamma^{i\star}_n, u^{i\star}_n,l^i\}}(\cdot).
\end{align}
Hence, we can also conclude that the empirical measure $\Lambda_n$ converges weakly through a subsequence to $\Lambda$, where the marginals on policies and actions coincide with the directing measure of the limiting infinitely exchangeable actions. 
\item [{\bf Step 2.}] {\bf Convergence of $N$-sensor expected exchangeable cost.}
Consider the optimal exchangeable subsequence of policies $P_{\pi}^{n\star}$ in Step 1 that converges to a limiting infinitely exchangeable policy $P_{\pi}^{\star}\in \LEX$ as $n\to \infty$.

Let 
\begin{align*}
    g(H, u^{i}; \gamma^{I_{i}^n\star}_n)&:=\int \nu(U^i=u^i|l^{i},P_{\pi}^{n\star},\mu^{n\star},I_{1}^n,\ldots, I_n^n)P(dl^i|H)\\
   &:=\int \gamma^{I_{i}^n\star}_n(u^i|l^i)P(dl^i|H).
\end{align*}
where $\mu^{n\star}$ is the empirical measure of random policies $\gamma^{1\star:n\star}_n$ induced by $P_{\pi}^{n\star}$.  Following step 1, we conclude that $g(H, u^{i}; \gamma^{I_{i}^n\star}_n)$ converges to 
\begin{align*}
     g(H, u^{i}; \gamma^{I_{i}\star})&:=\int \nu(U^i=u^i|l^{i},P_{\pi}^{\star},\mu^{\star},I_{1},I_2, \ldots)P(dl^i|H)\\
   &:=\int \gamma^{I_{i}\star}(u^i|l^i)P(dl^i|H)
\end{align*}
since $P_{\pi}^{n\star}$ and $\mu^{n\star}$ converges weakly to $P_{\pi}^{\star}$ and $\mu^{\star}$.

Following \cite[Lemma 1]{tsitsiklis1988decentralized} and \cite{shannon1967lower}, we have for any deterministic encoding policy $\gamma^{1:N}$ of the sensors,
\begin{align}\label{eq:bound-shannon}
   \inf_{\gamma^{0}} J^N_{\sf EE}(\gamma^{0:N})&\geq \min_{s\in [0,1]}\left[\frac{1}{N}\sum_{i=1}^{N}\log\left( \sum_{u^i\in \mathbb{U}} g(H_2, u^{i}; \gamma^{i})^{1-s} g(H_1, u^{i}; \gamma^{i})^{s} \right)\right] + \kappa_N,
   \end{align}
where 
\begin{align*}
    \kappa_N:=&- \frac{1}{N} \sqrt{2\sum_{i=1}^{N}\frac{d^2}{ds^2}\log\left( \sum_{u^i\in \mathbb{U}} g(H_2, u^{i}; \gamma^{i})^{1-s} g(H_1, u^{i}; \gamma^{i})^{s} \right)}\bigg{|}_{s^{\star}_N}+\frac{\log(\min\{P(H_1),P(H_2)\}/2)}{N}.
\end{align*}
We have
\begin{align*}
    &\limsup_{N\to \infty}\inf_{P_{\pi}\in \LEXN} J^N_{\sf EE}(\gamma^{0\star}, P_{\pi}^N)\\
    &\geq \lim_{n\to \infty} \mathbb{E}\left[\int \log\left( \sum_{u\in \mathbb{U}} g(H_2, u; \gamma)^{1-s_n^{\star}} g(H_1, u; \gamma)^{s_n^{\star}} \right) \Lambda_n(d\gamma,du,dl)\right]+\lim_{n\to \infty}\kappa_n\\ 
    &=\mathbb{E}\left[\int\log\left( \sum_{u\in \mathbb{U}} g(H_2, u; \gamma)^{1-s^{\star}} g(H_1, u; \gamma)^{s^{\star}} \right) \Lambda(d\gamma,du,dl)\right]\\
    &\geq \inf_{P_{\pi}^{R}\in \mathcal{P}(\Gamma^R)} \min_{s\in [0,1]} \int \log\left( \sum_{u\in \mathbb{U}} g(H_2, u; \gamma^{R})^{1-s} g(H_1, u; \gamma^{R})^{s} \right) P_{\pi}^{R}(d\gamma^R)\\
    &\geq \inf\limits_{P_{\pi} \in \LPRS}\limsup\limits_{N \to \infty} J^N_{\sf EE}(\gamma^{0\star}, P_{\pi,N})
    \end{align*}
    where the first inequality follows from Jensen's inequality, and $\kappa_n$ goes to zero since $$\sup_{\gamma^i,s}\frac{d^2}{ds^2}\log\left( \sum_{u^i\in \mathbb{U}} g(H_2, u^{i}; \gamma^{i})^{1-s} g(H_1, u^{i}; \gamma^{i})^{s} \right)<\infty$$ following \cite[Proposition 3]{tsitsiklis1988decentralized}.     
We have
\begin{flalign}
\inf\limits_{\gamma^0 \in \Gamma^{\sf FS}, P_{\pi} \in L}\limsup\limits_{N \to \infty} J^N_{\sf EE}(\gamma^0, P_{\pi,N})
&\geq\limsup\limits_{N \to \infty}\inf\limits_{\gamma^0 \in \Gamma^{\sf FS}_N, P_{\pi}^{N} \in \LCON}J^N_{\sf EE}(\gamma^0, P_{\pi}^N)\label{eq:info-red}\\
&=\limsup\limits_{N \to \infty}\inf\limits_{P_{\pi}^{N} \in \LEXN}J^N_{\sf EE}(\gamma^{0\star}, P_{\pi}^N)\label{eq:LN--LEXN}\\
&\geq \inf\limits_{P_{\pi} \in \LPRS}\limsup\limits_{N \to \infty} J^N_{\sf EE}(\gamma^{0\star}, P_{\pi,N})\label{eq:LPRS-Conv}\\
&\geq\inf\limits_{P_{\pi} \in L} \limsup\limits_{N \to \infty}  J^N_{\sf EE}(\gamma^{0\star}, P_{\pi,N})\label{eq:LPRS--L}\\
&\geq \inf\limits_{\gamma^0 \in \Gamma^{\sf FS}, P_{\pi} \in L}\limsup\limits_{N \to \infty} J^N_{\sf EE}(\gamma^0, P_{\pi,N}),
\end{flalign}
where \eqref{eq:info-red}  follows from exchanging limsup with infimum, and  Lemma \ref{Lem:info-optimal-N}. Equality \eqref{eq:LN--LEXN} follows from Lemma \ref{lemma:2-PI} which guarantees that the fusion center policy is a MAP rule $\gamma^{0\star}$ (under $I^{\sf FS}_N$ information structure), and Theorem \ref{the:exi-N}, which allows us to restrict the search for sensors' optimal policies to exchangeable ones. Inequality \eqref{eq:LPRS-Conv} follows from Step 2. Hence, the chain inequality above is an equality. 

\item [{\bf Step 3.}] {\bf  Optimality of infinitely exchangeable policies.}
By Theorem \ref{the:defin}, conditional on the directing measure $\mu^{\star}$, the limiting policy $P^{\star}_{\pi}$ corresponds to conditionally i.i.d. encoding policies $(\gamma^{1\star},\gamma^{2\star}, \ldots)$, i.e., $P^{\star}_{\pi}\in \LCOS$. Now, we show that $P^{\star}_{\pi}$ is optimal.

Let 
 \begin{align*}
     J^{\star}_{\mathrm{EE}}:=\inf_{P_\pi^R\in\mathcal P(\Gamma^R)}
\min_{s\in[0,1]}
\int_{\Gamma^R}
\log\left(
\sum_{u\in\mathbb U}
g(H_2,u;\gamma^R)^{1-s}
g(H_1,u;\gamma^R)^s
\right)P_\pi^R(d\gamma^R)
 \end{align*}

 Similarly to our analysis in Step 2, since the infinite exchangeable policy is conditionally i.i.d. with directing measure $\mu^\star$, we get
\begin{align*}
    \limsup_{N\to \infty}J^N_{\sf EE}(\gamma^0, P_{\pi,N}^{\star})&\geq J^{\star}_{\mathrm{EE}}
\end{align*}
 where $P_{\pi,N}^{\star}$ is the restriction of $P_{\pi}^{\star}$ to its first $N$ components. Additionally, following our analysis in Step 2, we have
\begin{align*}
    \mathbb{E}\left[\int\log\left( \sum_{u\in \mathbb{U}} g(H_2, u; \gamma)^{1-s^{\star}} g(H_1, u; \gamma)^{s^{\star}} \right) \Lambda(d\gamma,du,dl)\right]=J^{\star}_{\mathrm{EE}}
\end{align*}
 since for every \(\gamma\in\Gamma^R\) and \(s\in[0,1]\), 
$\log\left(
\sum_{u\in\mathbb U}
g(H_2,u;\gamma)^{1-s}
g(H_1,u;\gamma)^s
\right)\notag
\geq J^{\star}_{\mathrm{EE}} $, we have
 \begin{align}\label{eq:almost-sure-eq}
    \int\log\left( \sum_{u\in \mathbb{U}} g(H_2, u; \gamma)^{1-s^{\star}} g(H_1, u; \gamma)^{s^{\star}} \right) =J^{\star}_{\mathrm{EE}}
\end{align}
 $\Lambda$-almost surely.
 Hence, we have
\begin{align}
  &\limsup_{N\to \infty} \inf_{\gamma^{0}} J^N_{\sf EE}(\gamma^0,P_{\pi,N}^{\star})\notag\\&\leq \limsup_{N\to \infty}\frac{1}{N}\log\left(\mathbb{E} \left[\exp \left(N\int \log\left(  \sum_{u\in \mathbb{U}} g(H_2, u; \gamma)^{1-s^{\star}} g(H_1, u; \gamma)^{s^\star}\right)\Lambda(d\gamma,du,dl) \right)\right]\right)\label{eq:ineq-opt1}\\
   & \leq \limsup_{N\to \infty}\frac{1}{N}\log\left(\mathbb{E} \left[\exp \left(NJ^{\star}_{\mathrm{EE}} \right)\right]\right)\label{eq:ineq-opt2}\\
   &\leq J^{\star}_{\mathrm{EE}}\label{eq:ineq-opt3}\\
   &=\limsup_{N\to \infty}\inf_{P_{\pi}\in \LEXN} J^N_{\sf EE}(\gamma^{0\star}, P_{\pi}^N)\label{eq:ineq-opt4}
   \end{align}
where \eqref{eq:ineq-opt1} follows from the upperbound in \cite[Lemma 1]{tsitsiklis1988decentralized} 
for every \(N\ge 1\) and every deterministic encoding policy \(\gamma^{1:N}\),
\begin{align}
&\inf_{\gamma^0}J_{\sf EE}^N(\gamma^{0:N})\leq
\min_{s\in[0,1]}
\frac1N\sum_{i=1}^N
\log\left(
\sum_{u\in\mathbb U}
g(H_2,u;\gamma^i)^{1-s}
g(H_1,u;\gamma^i)^s
\right).
\end{align}
Inequality \eqref{eq:ineq-opt2} follows from \eqref{eq:almost-sure-eq} and \eqref{eq:ineq-opt4} follows from Step 2. This implies
\begin{align*}
    &\limsup_{N\to \infty}\inf_{P_{\pi}\in \LEXN} J^N_{\sf EE}(\gamma^{0\star}, P_{\pi}^N)= \limsup_{N\to \infty}J^N_{\sf EE}(\gamma^{0\star}, P_{\pi,N}^{\star})
\end{align*}
Hence, the limiting infinitely exchangeable policy $P^{\star}_{\pi}\in \LEX$ is asymptotically optimal. This completes the proof of part (i) since $\LCOS$ is a convex set with its extreme points lying in $\LPRS$ and the limiting expected cost of the representative problem is linear in randomized policies. Part (ii) follows from our derivation in part (i).
\end{itemize}

\subsection{Proof of Proposition \ref{the:thre-inf}}\label{APP:Prop1}
Part (i) directly follows from Theorem \ref{the:Exist-inf}. We now prove part (ii). 

Following Lemma \ref{Lem:info-optimal-N}, without any loss, we can restrict policies of the fusion center to $\Gamma^{0}_{\sf FS}$. Following Lemma \ref{lemma:1-EX}, we can restrict our search for optimal encoding policies to $\LEXN$ for $\mathcal{P}^{N}_{\sf E}$ without any loss of optimality. Let $P^{N\star}_{\pi} \in  \LEXN$ be an encoding optimal policy for $\mathcal{P}^{N}_{\sf E}$, i.e.,
\begin{align}\label{eq:asy32.1}
   \inf_{\gamma^{0}\in \Gamma^0_{\sf FS}} J^{N}_{\sf EE}(\gamma^{0},P^{ N\star}_{\pi}) \leq \inf\limits_{(\gamma^{0},P^{N}_{\pi}) \in \Gamma^0_{\sf FS} \times \LEXN}J^{N}_{\sf EE}(\gamma^{0},P^{ N}_{\pi})= J^{N\star}_{\sf EE}. 
\end{align}

 Let $P^{\star}_{\pi}\in \LPRS$ be an optimal encoding policy for $\mathcal{P}^{\infty}_{\sf E}$ with the restriction to the first $N$ components as $P^{
\star}_{\pi,N}$. Since $P_{\pi,N}^{\star}$ is symmetric and independent,   following from proof of Theorem \ref{the:Exist-inf}, we have
\begin{flalign}
\inf_{\gamma^0\in \Gamma^0_{\sf FS}} J^{N}_{\sf EE}(\gamma^0,P_{\pi,N}^{\star})\leq J^{N\star}_{\sf EE}+{\epsilon_{N}}
\end{flalign}
for some ${\epsilon_{N}}\geq 0$ such that ${\epsilon_{N}
}$ goes to $0$ as $N\to \infty$ which completes the proof.

\subsection{Proof of Theorem \ref{the:Robust}.}\label{app:rob}
{Since Assumption \ref{Assump:density} holds under $Q_n$ for all $n$, by Theorem \ref{the:thre}, an optimal solution exists. Let $(\gamma^{0\star}_n, \gamma^{1\star}_n, \ldots, \gamma^{N\star}_n)$ be an optimal solution under $Q_n$.  Let $\epsilon>0$ and choose $\gamma^{0\star:N\star}_{\epsilon}$ such that 
    \begin{align}
J^N(\gamma^{0\star:N\star}_{\epsilon};Q)-\inf_{\gamma^{0:N}} J^N(\gamma^{0:N};Q)\leq \frac{\epsilon}2 \label{eq:epsilon-opt}
    \end{align}
    We have
\begin{align}
&\left|\inf_{\gamma^{0:N}} J^N(\gamma^{0:N};Q_n)-\inf_{\gamma^{0:N}} J^N(\gamma^{0:N};Q)\right|\notag
    \\&= \left| J^N(\gamma^{0\star:N\star}_n;Q_n)-\inf_{\gamma^{0:N}} J^N(\gamma^{0:N};Q)\right|\notag\\
&\leq \left| J^N(\gamma^{0\star:N\star}_n;Q_n)-J^N(\gamma^{0\star:N\star}_{\epsilon};Q)\right|+\left|J^N(\gamma^{0\star:N\star}_{\epsilon};Q)-\inf_{\gamma^{0:N}} J^N(\gamma^{0:N};Q)\right|\label{eq:-Rob62}
    \end{align}
The second expression in \eqref{eq:-Rob62} is less than $\epsilon/2$ by \eqref{eq:epsilon-opt}. We now bound the first expression in \eqref{eq:-Rob62}. Suppose that $J^N(\gamma^{0\star:N\star}_n;Q_n)\geq J^N(\gamma^{0\star:N\star}_{\epsilon};Q)$, then
\begin{align}
J^N(\gamma^{0\star:N\star}_n;Q_n)-J^N(\gamma^{0\star:N\star}_{\epsilon};Q)\leq J^N(\gamma^{0\star:N\star}_{\epsilon};Q_n)-J^N(\gamma^{0\star:N\star}_{\epsilon};Q)\label{eq:-Rob62-i}
\end{align}
since $\gamma^{0\star:N\star}_n$ is optimal under $Q_n$, and hence, replacing it with $\gamma^{0\star:N\star}_{\epsilon}$ increases the expected cost and hence the difference. 
On the other hand, if $J^N(\gamma^{0\star:N\star}_n;Q_n)< J^N(\gamma^{0\star:N\star}_{\epsilon};Q)$, then
\begin{align}
J^N(\gamma^{0\star:N\star}_{\epsilon};Q)-J^N(\gamma^{0\star:N\star}_n;Q_n)\leq J^N(\gamma^{0\star:N\star}_n;Q)-J^N(\gamma^{0\star:N\star}_n;Q_n)+\frac{\epsilon}2\label{eq:-Rob62-ii}
\end{align}
Therefore, following \eqref{eq:-Rob62}--\eqref{eq:-Rob62-ii}, we get
    \begin{align*}
    &\left|\inf_{\gamma^{0:N}} J^N(\gamma^{0:N};Q_n)-\inf_{\gamma^{0:N}} J^N(\gamma^{0:N};Q)\right|\\
    &\leq \sup_{\gamma^{0:N}} \left|J^N(\gamma^{0:N};Q_n)- J^N(\gamma^{0:N};Q)\right| + \epsilon\\
    &\leq \sup_{\gamma^{0:N}} \bigg|  \int c(H,\gamma^{0}(u^{1:N})) \prod_{i=1}^{N}\gamma^{i}(du^{i}|y^i){Q}_n(dy^{1:N}|H_k)P(dH_k) - \\&\int c(H,\gamma^{0}(u^{1:N})) \prod_{i=1}^{N}\gamma^{i}(du^{i}|y^i){Q}(dy^{1:N}|H_k)P(dH_k)\bigg|+\epsilon\\
   &\leq\left\|PQ_n - PQ\right\|_{\sf TV} + \epsilon
\end{align*}
    This implies \eqref{eq:Continuity} since $\epsilon$ is arbitrary and $Q_n$ converges in total variation to $Q$ as $n\to \infty$.
    
Next, we establish robustness.  We have,
\begin{align*}
&\left|\ J^N(\gamma^{0\star:N\star}_n;Q)-\inf_{\gamma^{0:N}} J^N(\gamma^{0:N};Q)\right|\\
&\leq \left|\ J^N(\gamma^{0\star:N\star}_n;Q)-J^N(\gamma^{0\star:N\star}_n;Q_n)\right|+\left|\ J^N(\gamma^{0\star:N\star}_n;Q_n)-\inf_{\gamma^{0:N}} J^N(\gamma^{0:N};Q)\right|
\end{align*}
The second term on the right-hand side goes to zero by the established continuity; we now show that the first one also goes to zero as $n\to \infty$. We have
\begin{align*}
&\left|\ J^N(\gamma^{0\star:N\star}_n;Q)-J^N(\gamma^{0\star:N\star}_n;Q_n)\right|\\
&\leq \sup_{\gamma^{0:N}} \bigg|  \int c(H,\gamma^{0}(u^{1:N})) \prod_{i=1}^{N}\gamma^{i}(du^{i}|y^i){Q}_n(dy^{1:N}|H_k)P(dH_k) - \\
&\int c(H,\gamma^{0}(u^{1:N})) \prod_{i=1}^{N}\gamma^{i}(du^{i}|y^i){Q}(dy^{1:N}|H_k)P(dH_k)\bigg|\\
   &\leq\left\|PQ_n - PQ\right\|_{\sf TV}\to 0
\end{align*}
which completes the proof.

\subsection{Proof of Theorem \ref{the:Approx-rep}.}\label{APP:approx-proof}

We first prove a lemma on the joint continuity of the cost in policies and the parameter $s$, which we use to prove Theorem \ref{the:Approx-rep}.
\begin{lemma}\label{lem:cont-rep}
Let $P_i(dl):=P(dl\mid H_i)$ for $i=1,2$. Define $P:=(P_1+P_2)/2$.
Let \(\Gamma^R\) be such that $\gamma(dl,du)=P(dl)\gamma_\pi(du\mid l)$, endowed with the Young topology. Under Assumptions \ref{Assump:density} and \ref{Assump:inf-g}, the map
\[
(\gamma,s)
\mapsto
\log\left(
\sum_{u\in\mathbb U}
g(H_2,u;\gamma)^{1-s}
g(H_1,u;\gamma)^s
\right)
\]
is jointly continuous on \(\Gamma^R\times[0,1]\).
\end{lemma}

\begin{proof}[Proof of Lemma \ref{lem:cont-rep}]
The square-integrability condition in Assumption \ref{Assump:inf-g} implies that
$0<l^R<\infty$ $P$-a.s. Now let
$h_i(l):=\frac{dP_i}{dP}(l)$ for 
$i=1,2$. For every \(u\in\mathbb U\),
\[
\begin{aligned}
g(H_i,u;\gamma^R)
&=
\int_{\mathbb L}\gamma_{\pi}^R(u\mid l)P_i(dl)  \\
&=
\int_{\mathbb L\times\mathbb U}
h_i(l)\mathbb{I}_{\{u\}}(a)\gamma^R(dl,da).
\end{aligned}
\]
Suppose that $\gamma_n^R\to\gamma^R$
in the Young topology. Since \(h_i\) is bounded and measurable in \(l\), and
since \(\mathbb U\) is finite, the map
$a\mapsto \mathbb{I}_{\{u\}}(a)$ is continuous on \(\mathbb U\). Therefore, by the definition of convergence in the Young topology,
$g(H_i,u;\gamma_n^R)\to g(H_i,u;\gamma^R)$ for $i=1,2$ and $u\in\mathbb U$. Thus the finite-dimensional vectors
$\bigl(g(H_i,u;\gamma_n^R)\bigr)_{i=1,2,\ u\in\mathbb U}$
converge to $\bigl(g(H_i,u;\gamma^R)\bigr)_{i=1,2,\ u\in\mathbb U}$.

Also, if $g(H_1,u;\gamma^R)=0$,
then $\int_{\mathbb L}\gamma^R_{\pi}(u\mid l)P_1(dl)=0$.
This implies $\gamma^R_{\pi}(u\mid l)=0$
$P_1$-a.s. Since $P_2$ is absolutely continuous with respect to $P_1$ by Assumptions \ref{Assump:density} and \ref{Assump:inf-g}, we also have $\gamma^R_{\pi}(u\mid l)=0$
$P_2$-a.s., and hence
$g(H_2,u;\gamma^R)=0$.
The reverse implication follows similarly. Therefore,
\[
g(H_1,u;\gamma^R)=0
\quad\Longleftrightarrow\quad
g(H_2,u;\gamma^R)=0,
\qquad u\in\mathbb U.
\]
The same argument applies to every \(\gamma_n^R\).

Now let \(s_n\to s\) in \([0,1]\) and denote
$p_{n,u}:=g(H_1,u;\gamma_n^R)$, $q_{n,u}:=g(H_2,u;\gamma_n^R)$, $p_{u}:=g(H_1,u;\gamma^R)$, and
$q_{u}:=g(H_2,u;\gamma^R)$.
We have \(p_{n,u}\to p_u\) and \(q_{n,u}\to q_u\) for every \(u\in\mathbb U\). If \(s\in(0,1)\), then
$(q,p,s)\mapsto q^{1-s}p^s$
is continuous on \([0,1]^2\times(0,1)\). Hence
$q_{n,u}^{1-s_n}p_{n,u}^{s_n}
\to
q_u^{1-s}p_u^s$
for every \(u\in\mathbb U\). Since \(\mathbb U\) is finite,
\[
\sum_{u\in\mathbb U}
q_{n,u}^{1-s_n}p_{n,u}^{s_n}
\to
\sum_{u\in\mathbb U}
q_u^{1-s}p_u^s.
\]
It remains to check the endpoints. Suppose \(s=0\). If \(q_u>0\) for $u\in \mathbb{U}$, then \(p_u>0\), and therefore
\[
q_{n,u}^{1-s_n}p_{n,u}^{s_n}\to q_u.
\]
If \(q_u=0\) for $u\in \mathbb{U}$, then also \(p_u=0\). Since \(q_{n,u}\to0\) and \(p_{n,u}\to0\),
and since \(s_n\to0\), for all sufficiently large \(n\) we have $1-s_n\ge 1/2.$ Since \(0\le p_{n,u}\leq 1\), so $0\le p_{n,u}^{s_n}\leq 1$.
Therefore,
\[
0\le q_{n,u}^{1-s_n}p_{n,u}^{s_n}
\le q_{n,u}^{1/2}\to0.
\]
Therefore,
\[
\sum_{u\in\mathbb U}
q_{n,u}^{1-s_n}p_{n,u}^{s_n}
\to
\sum_{u\in\mathbb U}q_u
=
1.
\]
The endpoint \(s=1\) is analogous. Therefore, for every \(s\in[0,1]\),
\[
\sum_{u\in\mathbb U}
g(H_2,u;\gamma_n^R)^{1-s_n}
g(H_1,u;\gamma_n^R)^{s_n}
\to
\sum_{u\in\mathbb U}
g(H_2,u;\gamma^R)^{1-s}
g(H_1,u;\gamma^R)^s.
\]
The limiting sum is strictly positive, and hence composition with the logarithm is continuous. Therefore,
\[
\log\left(
\sum_{u\in\mathbb U}
g(H_2,u;\gamma_n^R)^{1-s_n}
g(H_1,u;\gamma_n^R)^{s_n}
\right)
\to
\log\left(
\sum_{u\in\mathbb U}
g(H_2,u;\gamma^R)^{1-s}
g(H_1,u;\gamma^R)^s
\right).
\]
This completes the proof of the lemma.
\end{proof}

\begin{proof}[Proof of Theorem \ref{the:Approx-rep}]
Let 
\begin{align*}
J_{\sf EE}^{\star,n}
:=
\inf_{\phi_n^R\in\Phi_n^R}
\min_{s\in [0,1]}
\log\left(
\sum_{u\in\mathbb U}
g(H_2,u;\phi_n^R)^{1-s}
g(H_1,u;\phi_n^R)^s
\right).
\end{align*}
Define, for
$g=\bigl(g(H_i,u)\bigr)_{i=1,2,\ u\in\mathbb U}\in\mathcal G$,
the function
\[
F(g)
:=
\min_{s\in[0,1]}
\log\left(
\sum_{u\in\mathbb U}
g(H_2,u)^{1-s}g(H_1,u)^s
\right).
\]
Following the denseness argument in \cite[Proof of Theorem 5.1]{yuksel2023borkar} under the Young topology for the representative encoder policy, we get that for every representative encoder policy
\(\gamma^R\in\Gamma^R\), there exists a sequence of finite-cell deterministic
encoders \(\phi_n^R\in\Phi_n^R\) that converges to it in the Young topology, and hence, following the same reasoning as that used in the proof of Lemma \ref{lem:cont-rep}, we get
\begin{align}\label{eq:g-conv}
\bigl(g(H_i,u;\phi_n^R)\bigr)_{i=1,2,\ u\in\mathbb U}
\to
\bigl(g(H_i,u;\gamma^R)\bigr)_{i=1,2,\ u\in\mathbb U}.
\end{align}
Following the proof of Lemma \ref{lem:cont-rep}, for each fixed \(s\in[0,1]\), the map
\[ g \mapsto \sum_{u\in\mathbb U} g(H_2,u)^{1-s}g(H_1,u)^s \] is well-defined and continuous, and hence, the minimum is attainable. Thus, following Lemma \ref{lem:cont-rep}, $F$ is continuous since, for a general function $\hat{F}(x,y)$ which is jointly continuous, $\min_{x \in A}\hat{F}(x,y) =:G(y)$ is continuous. Following Theorem \ref{the:Exist-inf}, an optimal encoder policy exists. Let $g^{\star}\in \mathcal{G}$ correspond to the optimal encoder policy. Then, an analogous argument as that used for \eqref{eq:g-conv} implies that there exists $g^{\star}_n\in \mathcal{G}_n$ that converges to $g^{\star}$ as $n\to \infty$. This implies that $F(g^{\star}_n)$ converges to 
$F(g^{\star})$ as $n\to \infty$. Hence,
\begin{align*}
    \limsup_{n\to \infty}J^{\star,n}_{\sf EE}&=\limsup_{n\to \infty}\inf_{g_n \in \mathcal{G}_n}F(g_n)\\
    &\leq \lim_{n\to \infty}F(g^{\star}_n)\\
    &=F(g^{\star})\\
    &=J^{\star}_{\sf EE},
\end{align*}
where the first inequality follows from the fact that $g^{\star}_n$ is possibly suboptimal.
On the other hand, we have $\mathcal{G}_n\subset \mathcal{G}$ since \(\phi_n^R\) is constant on every quantization cell \(A_{n,k}\), and hence every vector in \(\mathcal{G}_n\) is also generated by an admissible encoder under the original measure \(P\). Hence, we conclude that $J^{\star,n}_{\sf EE}\geq J^{\star}_{\sf EE}$ for all $n$, and hence,
\begin{align*}
    \liminf_{n\to \infty}J^{\star,n}_{\sf EE}\geq J^{\star}_{\sf EE}.
\end{align*}
Thus,
\begin{align*}
    J^{\star}_{\sf EE}\leq  \liminf_{n\to \infty}J^{\star,n}_{\sf EE} \leq \limsup_{n\to \infty}J^{\star,n}_{\sf EE}\leq J^{\star}_{\sf EE}.
\end{align*}
Hence,
\begin{align}\lim_{n\to \infty}J^{\star,n}_{\sf EE}=J^{\star}_{\sf EE}.\label{eq:converg-n}
\end{align}

For the finite grid \(S_m\subset[0,1]\), define
\[
F_m(g)
:=
\min_{s\in S_m}
\log\left(
\sum_{u\in\mathbb U}
g(H_2,u)^{1-s}g(H_1,u)^s
\right).
\]
We have $F(g)\le F_m(g)$ for all
$g\in\mathcal G_n$, and hence, $J_{\sf EE}^{\star,n}
\le
J_{\sf EE}^{\star,n,m}$. Since the map
\[
(g,s)\mapsto
\log\left(
\sum_{u\in\mathbb U}
g(H_2,u)^{1-s}g(H_1,u)^s
\right)
\]
is continuous on the compact set \(\mathcal{G}_n\times[0,1]\), it is uniformly continuous.
As the mesh of \(S_m\) tends to zero, it follows that
\[
\sup_{g\in \mathcal{G}_n}|F_m(g)-F(g)|\to0.
\]
Therefore, for every fixed \(n\),
$\lim_{m\to\infty}J_{\sf EE}^{\star,n,m}
=
J_{\sf EE}^{n\star}$, which completes the proof by \eqref{eq:converg-n}.
\end{proof}
}

\subsection{Proof of Theorem \ref{the:Applic}}\label{APP:Thm6}
{
Part (i) follows from an analogous argument used in the proofs of Lemma \ref{lemma:1-EX} and Theorem \ref{the:thre}. We first show that the space of policies is compact under total variation and that the expected throughput is continuous; hence, we conclude that an optimal policy exists.  

For large $t$, by our assumption, $L_t^{i}> 0$ for all $i\in \mathbb{N}$. Hence, we can define the randomized transmission policies $\lambda^N_i$ as a probability measure on the action space $\{0,1\}$ for each user and a finite nonnegative measure taking values in $[0,\frac{\bar{\alpha}}{N}]$ for large $N$. As a result, the space of randomized transmission policies is compact under total variation for each user. Following an analogous argument to step 2 of the proof of Lemma \ref{lem:compactness of TS}, we can show that the space of randomized transmission policies among all users is also compact under total variation. We next show that the expected throughput is continuous in randomized policies. Suppose $\lambda^{N,k}_{i}$ converges to $\lambda^{N}_{i}$ as $k\to \infty$ in total variation for all $i=1, \ldots, N$. Then, for large $t$, by our assumption, $L_t^{i}> 0$ for all $i\in \mathbb{N}$, and hence, 
    \begin{align*}
    \left|J^N_t(\pi^N_k)-J^N_t(\pi^N)\right|
    &=
    \bigg|\int \sum_{i=1}^N 1_{\{u^i_t=1\}}
\prod_{j\neq i}(1-1_{\{u^{j}_t=1\}})\prod_{i=1}^{N}\lambda_i^{N,k}(du_t^{i})\\
    &-\int \sum_{i=1}^N 1_{\{u^i_t=1\}}
\prod_{j\neq i}(1-1_{\{u^{j}_t=1\}})\prod_{i=1}^{N}\lambda_i^{N}(du_t^{i})\bigg|\\
&\leq \left\|\prod_{i=1}^{N}\lambda_i^{N,k}(du_t^{i})-\prod_{i=1}^{N}\lambda_i^{N}(du_t^{i}) \right\|_{\sf TV}
    \end{align*}
which goes to zero by an analogous argument to step 2 of the proof of Lemma \ref{lem:compactness of TS}. The rest of the proof follows from a similar argument as that used in the proof of Lemma \ref{lemma:1-EX}, by constructing an exchangeable policy as an average of all its permutations, since for finite $N$ the throughput and dynamics are exchangeable. This implies that there exists an optimal policy ${\pi}^{N\star}$ that is exchangeable, concluding the proof of part (i).}

{The proof of Part (ii) proceeds in two steps:
\begin{itemize}
  \item [{\bf Step 1.}] {\bf Asymptotic upperbound of the optimal expected reward.}
 Following part (i), for finite but large $N$, there exists an optimal policy ${\pi}^{N\star}$ that is exchangeable. By compactness of $[0,\bar{\alpha}]$, and analogous argument as used in the proof of Theorem \ref{the:Exist-inf}, there exists a subsequence of optimal finite-$N$ policies $\{{\pi}^{n\star}\}_n$ that converges to an infinitely exchangeable policy $\pi^{\star}$ as $n\to \infty$. We note that the policy $\pi^{\star}$ exists as a probability measure on $\{\alpha^i\}_{i=1}^{\infty}$, but since $N\to \infty$, the low transmission policy automatically leads to no transmission policy due to the boundedness assumption on $\alpha$.

For every $t\geq 0$, define
 \[
T^N_t(\pi^{N\star}):=\mathbb{E}^{\pi^{N\star}}\!\left[\sum_{i=1}^N 1_{\{U^i_t=1\}} \prod_{j\neq i} (1-  1_{\{U^j_t=1\}})\right]
\]
where in the above, we assume that $L_t^{i}> 0$ for all $i\in \mathbb{N}$. 
For every $t\geq 0$ and $N\in \mathbb{N}$, we have
\begin{align*}
J^N_t(\pi^{N\star})-T_t^N(\pi^{N\star}) \leq 0.
\end{align*}
 
Let $\pi|_{N}$ be restriction of $\pi$ on the first $N$ components. Hence,
\begin{align*}
\liminf_{N\to \infty} \sup_{\pi^N} J^N_t(\pi^{N})
&=\liminf_{N\to \infty} J^N_t(\pi^{N\star})\\
&\leq \liminf_{N\to \infty} T^N_t(\pi^{N\star})
\end{align*}
Under $\pi^{N\star}$, we get
\begin{align*}
\widehat{T}^N(\pi^{N\star})&:=T^N_t(\pi^{N\star})\\
&=\mathbb{E}^{\pi^{N\star}}\left[\!\mathbb{E}\left[\sum_{i=1}^N 1_{\{U^i_t=1\}} \prod_{j\neq i} (1-  1_{\{U^j_t=1\}})\bigg{|} \alpha_1^N, \ldots, \alpha_N^N\right]\right]\\
&=\mathbb{E}^{\pi^{N\star}}\left[\sum_{i=1}^N \frac{\alpha_i^N}{N}\prod_{j\neq i}\left(1-\frac{\alpha_j^N}{N}\right)\right].
\end{align*}
Let $\mu_N$ be the empirical measures on $\alpha$s under optimal policy $\pi^{N\star}$, i.e., 
\[\mu_N(\cdot)=\frac{1}{N}\sum_{i=1}^{N} \delta_{\alpha_i^{N\star}}(\cdot)\]
and 
\[
m_N:=\frac{1}{N}\sum_{i=1}^N \alpha_i^{N\star}=\int \alpha \mu_N(d\alpha).
\] 

Using $1-x\le e^{-x}$, we obtain
\[
\prod_{j\neq i}\left(1-\frac{\alpha_j^N}{N}\right)
\le
\exp\!\left(-\frac{1}{N}\sum_{j\neq i}\alpha_j^N\right)
=
\exp\!\left(-m_N+\frac{\alpha_i^N}{N}\right).
\]
Let $$g^N(\alpha)=\alpha \exp\!\left(\frac{\alpha}{N}\right).$$
Since $\alpha$ takes value from a bounded set $[0, \bar{\alpha}]$, $g^N(\alpha)\leq \alpha \exp\!\left(\frac{\bar{\alpha}}{N}\right)$. Then, we have
\begin{align*}
\liminf_{N\to \infty}\widehat{T}^N(\pi^{N\star})
&\le
\liminf_{N\to \infty}\mathbb{E}^{\pi^{N\star}}\left[\frac{1}{N}\sum_{i=1}^N \alpha_i^N
\exp\!\left(-m_N+\frac{\alpha_i^N}{N}\right)\right]
\\
&=\liminf_{N\to \infty}\mathbb{E}^{\pi^{N\star}}\left[e^{-m_N}\frac{1}{N}\sum_{i=1}^N \alpha_i^N
\exp\!\left(\frac{\alpha_i^N}{N}\right)\right]\\
&=\liminf_{N\to \infty}\mathbb{E}^{\pi^{N\star}}\left[e^{-m_N}\frac{1}{N}\sum_{i=1}^N g^N(\alpha_i^N)
\right]\\
&=\liminf_{N\to \infty}\mathbb{E}^{\pi^{N\star}}\left[e^{-m_N}\int g^N(\alpha) \mu_N(d\alpha)
\right]\\
&\leq\liminf_{N\to \infty}\exp\!\left(\frac{\bar{\alpha}}{N}\right)\mathbb{E}\left[e^{-m_N}\int \alpha\mu_N(d\alpha)
\right]\\
&= \liminf_{N\to \infty}\mathbb{E}\left[m_Ne^{-m_N}\right].
\end{align*}
Recall that from part (i), for every finite $N$, the optimal policy $\pi^{N\star}$ is exchangeable with the induced empirical measure $\mu_N$, which is a random probability measure. Hence, the empirical average $m_N$ is a $[0, \bar{\alpha}]$-valued random variable. Since by compactness, a subsequence of exchangeable random variables $\{\alpha^{n\star}_{i}\}_{i=1}^{n}$ converges in distribution to an infinitely exchangeable random variables $\{\alpha^{\star}_{i}\}_{i=1}^{\infty}$ (i.e., $(\alpha_{1}^{n\star},\dots, \alpha_{m}^{n\star})$ converges in distribution to $(\alpha_{1}^{\star},\dots, \alpha_{m}^{\star})$ for each $m\geq 1$ as $n\to \infty$), following \cite[Proposition 7.20]{aldous2006ecole}, the subsequence of their empirical measures $\{\mu_n\}_n$ converges weakly to a random measure $\mu$, the directing measure of the infinite exchangeable limiting random variables $\{\alpha^{\star}_{i}\}_{i=1}^{\infty}$. This yields that the subsequence of empirical averages $\{m_n\}_{n}$ of $\{\alpha^{n\star}_{i}\}_{i=1}^{\infty}$ converges in distribution to a random variable 
\[m=\int \alpha \mu(d\alpha)\]
as $n\to \infty$. Hence,
\begin{align*}
   \liminf_{N\to \infty}\mathbb{E}\left[m_Ne^{-m_N}\right]&\leq \lim_{n\to \infty}\int m_ne^{-m_n}\mathbb{P}(d\omega)\\
&= \mathbb{E}\left[ me^{-m}\right],
\end{align*}
where the last equality follows from dominated convergence theorem since $m_ne^{-m_n}$ is bounded by $e^{-1}$ following the fact that $f(x)=x e^{-x}$ has a unique maximizer in $[0,\infty)$ at $x=1$.

\item [{\bf Step 2.}] {\bf Optimality of infinitely exchangeable policies.} Suppose that $\pi^{\star}$ is infinitely exchangeable. By de Finetti's theorem, this policy $\pi^{\star}$ has a directing measure $\mu$ such that, conditional on $\mu$, the users' transmission policies are i.i.d.\ with attempt probability $\alpha/N$ with $\alpha$ distributed as $\mu$. Suppose that asymptotically as $t\to \infty$, each user has a packet to send.

Hence, for large $t$, conditional on $\mu$, the one-slot throughput is
\begin{align*}
J^N_t(\pi^{\star})&=\mathbb{E}\left[\mathbb{E}\left[\sum_{i=1}^N \frac{\alpha^i}{N}\prod_{j\not= i}\left(1-\frac{\alpha^j}{N}\right) \bigg| \mu\right]\right]\\
&=\mathbb{E}\left[\sum_{i=1}^N \frac{m}{N}\left(1-\frac{m}{N}\right)^{N-1} \right]
\\
&=
\mathbb{E}\left[ m\left(1-\frac{m}{N}\right)^{N-1} \right].
\end{align*}
where $m=\int \alpha \mu(d\alpha)$.
Since $0\le m\le \bar\alpha$, for $N>\bar\alpha$, we have
\[
0 \le m\left(1-\frac{m}{N}\right)^{N-1}\le m\le \bar\alpha,
\]
so by the dominated convergence theorem
\[
\lim_{N\to\infty} J^N_t(\pi^{\star})
=
\mathbb E\!\left[\lim_{N\to\infty} m\left(1-\frac{m}{N}\right)^{N-1}\right]
=
\mathbb E[m e^{-m}]\leq e^{-1}.
\]

The function
\[
f(\alpha):=\alpha e^{-\alpha}, \qquad \alpha\ge 0,
\]
is maximized at $\alpha=1$, with
\[
\sup_{\alpha\ge 0} \alpha e^{-\alpha}=e^{-1}.
\]
Hence,
\[
\mathbb E[m e^{-m}]
\le
\sup_{\alpha\ge 0} \alpha e^{-\alpha}
=
e^{-1}.
\]

Since $f(\alpha)=\alpha e^{-\alpha}$ has a unique maximizer at $\alpha=1$, the equality
\[
\mathbb E[m e^{-m}] = e^{-1}
\]
holds if and only if
\[
m=1 \qquad \text{almost surely}.
\]
In particular, when the directing measure is degenerate at $1$, the bound is achieved, i.e., $\mu(\cdot)=\delta_{1}(\cdot)$. Hence, the near-optimal policy is symmetric and independent, with each user having a per-user transmission probability of $1/N$.
    \end{itemize}
}

\bibliographystyle{IEEEtran}

\end{document}